\documentclass[10pt]{article}
\usepackage{amsmath}
\usepackage{graphicx}
\usepackage{amsfonts}
\usepackage{amssymb}
\usepackage{times}

\begin{document}

\begin{center}
.

\vspace{0.5in}

{\LARGE \textbf{A Coaxial Vortex Ring Model for Vortex Breakdown}}

\vspace{0.3in}

\textsf{Denis Blackmore}

\textsf{Department of Mathematical Sciences and}

\textsf{Center for Applied Mathematics and Statistics}

\textsf{New Jersey Institute of Technology}

\textsf{Newark, NJ 07102-1982}

\textsf{deblac@m.njit.edu }

$\ast\;\ast\;\ast$\textsf{\ }

\textsf{Morten Br\o ns}

\textsf{Department of Mathematics}

\textsf{Technical University of Denmark}

\textsf{DK-2800 Kgs. Lyngby, Denmark}

\smallskip\textsf{m.brons@mat.dtu.dk}

$\ast\;\ast\;\ast$\textsf{\ }

\textsf{Arnaud Goullet }

\textsf{Department of Mathematical Sciences}

\textsf{New Jersey Institute of Technology}

\textsf{Newark, NJ 07102-1982}

\textsf{abg3@njit.edu }
\end{center}

\vspace{0.2in}

\noindent\textbf{ABSTRACT: }A simple - yet plausible - model for B-type
vortex breakdown flows is postulated; one that is based on the immersion of
a pair of slender coaxial vortex rings in a swirling flow of an ideal fluid
rotating around the axis of symmetry of the rings. It is shown that this
model exhibits in the advection of passive fluid particles (kinematics) just
about all of the characteristics that have been observed in what is now a
substantial body of published research on the phenomenon of vortex
breakdown. Moreover, it is demonstrated how the very nature of the fluid
dynamics in axisymmetric breakdown flows can be predicted and controlled by
the choice of the initial ring configurations and their vortex strengths.
The dynamic intricacies produced by the two ring + swirl model are
illustrated with several numerical experiments.

\bigskip


\medskip


\bigskip

\section{introduction}

Vortex breakdown is officially five decades old this year: Although there is
evidence to suggest that this phenomenon was observed several centuries ago
in various atmospheric contexts, the first officially recorded description
appears to be that in the experimental paper of Peckham \& Atkinson \cite%
{PeckAt}. This was followed hard upon by the more vortex breakdown focused
investigations of Elle \cite{elle} and Lambourne \& Bryer \cite{LB}, and
what appears to be the first real attempt at analyzing the phenomenon in
Squire \cite{squire}. And even after half a century of research, interest in
vortex breakdown has not diminished; in fact, it has probably increased
significantly over the last fifteen or so years.

The early work on vortex breakdown quickly captured the attention of a
substantial group of talented fluid mechanicians attracted by the
importance, and challenge of unraveling the secrets of this intriguing and
rather mysterious phenomenon. With this intense level of scrutiny, it did
not take long to realize from the experimental evidence that there were
apparently two distinct types of vortex breakdown configurations: an
ostensibly axisymmetric \emph{bubble type} (\emph{B-type}) manifestation;
and a definitely asymmetric \emph{spiral type} (\emph{S-type}) structure
often observed as emanating from the trailing edge of a B-type breakdown
form. Moreover, inspired by Squire's attempt to formulate a theoretical
explanation of the phenomena in terms of variationally based critical state
criteria, many similar and quite different theoretical descriptions were
propounded. The three main types of theoretical characterizations employed
critical state concepts, boundary layer separation analogies, and long wave
hydrodynamic instability formulations, and a variety of flow properties and
characteristic combinations of parameters, such as the Rossby number, were
identified as playing key roles in the formation of vortex breakdown
configurations. Outstanding examples of these contributions, in addition to
Squire's pioneering analysis, can be found in the work of Benjamin \cite%
{benj}, Ludwieg \cite{ludw}, Hall \cite{hall}, Leibovich \cite{leib},
Escudier \& Keller \cite{esckell}, Trigub \cite{trigub}, Spall \emph{et al}.
\cite{SGG}, Berger \cite{berg}, Rusak \emph{et al}. \cite{RWW}, Gelfgat
\emph{et al.} \cite{GBS}, and Krause \cite{krause2, krause3}. These and
other theoretical investigations have done a great deal to expand our
understanding of vortex breakdown, especially the B-type variety.
Nevertheless, the venerable enigma that is vortex breakdown has yet to yield
to a complete and universally accepted exposition.

Soon after the quest for a theoretical explanation of vortex breakdown
began, a host of efforts directed at identifying vortex breakdown flows via
computational techniques applied to both the Navier-Stokes and Euler
equations were undertaken. Early results in this vein were somewhat
inconclusive, but improvements in numerical methods and major advances in
computers have produced results that are quite compelling. Notable examples
of this important aspect of vortex breakdown research can be found in
Grabowski \& Berger \cite{graberg}, Krause \emph{et al}. \cite{KSH}, Hafez
\emph{et al}. \cite{HKS}, Spall \& Gatski \cite{SG}, Breuer \cite{BREU},
Weimer \cite{WEIMER}, and Krause \cite{krause1}. All of these numerical
investigations have benefitted from and been complemented by a number of
quite sophisticated experimental studies of vortex breakdown that began in
earnest around the same time, with Sarpkaya \cite{sarp1}, Faler \& Leibovich
\cite{FL}, and Uchida \emph{et al}. \cite{UNO} among the best examples of
the earlier experimentally based investigations.

Holmes \cite{holmes} was instrumental in starting a relatively new trend in
research in fluid mechanics in general and vortex breakdown flows in
particular: During the last twenty years, many investigations of a
theoretical, computational, or experimental nature have been conducted from
the perspective of modern dynamical systems theory. Terminology such as
structural stability, stable and unstable manifolds, homoclinic and
heteroclinic orbits, and strange attractors have now become the lingua
franca for a large segment of vortex breakdown investigations. Prime
examples of this trend can be found in the work of Blackmore \cite{dbVB},
Blackmore \& Knio \cite{dbokVB}, Br\o ns and his collaborators \cite{BVS,
BB, BSSZ}, Gelfgat \emph{et al}. \cite{GBS}, Serre \& Bontoux \cite{serbon},
S\o rensen \& Christiansen \cite{SC}, and Sotiropoulos and his
coinvestigators \cite{sotir1, sotir2}. Viewing the results of this dynamical
systems approach, it is safe to say that it has already proved quite useful
in apprehending some of the more complex features of vortex breakdown flows,
and has great potential for catalyzing future leaps of understanding in this
area.

Dynamical systems theory certainly inspired the work presented here on our
vortex breakdown model, which was first adumbrated in \cite{dbokVB}. But the
idea of using a pair of circular rings was also influenced and guided by the
experimental and computational literature and interaction with some of the
leading practitioners in these areas. On the one hand, we knew from the
dynamical systems perspective that, among other things, two rings are
capable of generating the kinds of chaotic kinematics that appear to be
quite common in vortex breakdown flows \cite{bb, blkn1, blkn2, BWC, HSC},
but on the other hand, our choice was informed by the apparent ubiquity of
two coaxial vortex rings in numerical and experimental investigations.
Naturally, our decision to concentrate on an ideal fluid based model was
motivated by a desire to keep the model simple and also have the consequent
rich symplectic structure at our disposal; the advantages of which are, we
believe, amply demonstrated in the sequel.

The development and analysis of our two coaxial vortex ring model for B-type
vortex breakdown unfolds in this paper as follows: In Section 2, we
formulate the mathematical model for the motion of a pair of coaxial vortex
rings immersed in a swirling ideal fluid flow along the axis of symmetry of
the rings in the context of Hamiltonian dynamics, wherein we consider such
features as the complete integrability of the resulting system. Next, in
Section 3, we derive the Hamiltonian dynamical equations for the motion
induced in passive fluid particles (kinematics) by the dynamics of the rings
in the swirling flow. The associated kinematic equations are formulated in
two ways: in terms of a passive third ring with zero vortex strength; and
directly as a (time-dependent) Hamiltonian one-degree-of-freedom system. It
is noted that the resulting kinematics is not completely integrable, so it
is capable of exhibiting chaotic flow regimes. In addition, we introduce the
approach to studying the dynamics via a Poincar\'{e} map associated to
successive intersections of the streamlines with a fixed meridian
half-plane. Then the usual aspects of the dynamical system for the ring
motion, including classification of stationary points, are analyzed in
Section 4. This is followed in Section 5 with an analogous investigation of
the kinematic equations, which incorporates a brief analysis of homoclinic
and heteroclinic orbits connecting the several fixed points.

Since it was found that the kinematics for the case of fixed rings is too
regular to subsume the kinds of exotic streamline patterns that one expects
in vortex breakdown on the basis of several careful experimental and
numerical investigations, we formulate, in Section 6, the model equations in
terms of small (quasiperiodic but not necessarily periodic) oscillations of
the rings. Particular attention is paid to obtaining asymptotic expansions
of the heteroclinic orbit comprising the outer boundary of the bubble for
the trace of the kinematics in a meridian half-plane. Then, in Section 7, we
prove using Melnikov's method that small oscillations in the rings produce
chaotic streamline configurations of the type observed - but not proved - in
earlier studies of vortex breakdown phenomena. This is followed by Section 8
in which a variety of numerical experiments are run to illustrate the
dependence of the kinematics on various model parameters such as the swirl
strength and relative vortex ring strengths. The nature of transitions to
chaos as the size of the ring oscillations is increased is given special
emphasis. Finally, we conclude in Section 9 with a discussion of our model,
its plausibility, and its efficacy in predicting and controlling vortex
breakdown behavior. We also indicate some natural further research
directions suggested by the results obtained in our investigation.

\section{The Two Ring Model}

First, we consider the dynamics of two coaxial vortex moving in an ideal (=
inviscid and incompressible) fluid initially at rest in $\mathbb{R}^{3}$.
The axis of symmetry of the vortex rings is chosen to be the $x$-axis, and
their respective positive strengths are denoted as $\Gamma_{1}$ and $%
\Gamma_{2}$. Owing to the axisymmetry of the motion of these rings, their
positions are completely determined by their intersection points with any
meridian half-plane containing the $x$-axis. In order to fix ideas, we
choose this half-plane, which we denote by $\mathfrak{H}$, to be the $xy$%
-plane with $y\geq0$ ($y\geq0,z=0$), so that the motion of the rings is
characterized by their respective points of intersection $%
(x_{1},y_{1}):=(x_{1},y_{1},0)$ and $(x_{2},y_{2}):=(x_{2},y_{2},0)$ with
the designated half-plane.

Defining $s:=r^{2}:=y^{2}+z^{2}$, which reduces to $y^{2}$ on $\mathfrak{H}$%
, it is well known \cite{bb, blkn1, blkn2, BWC, BVS, BB} that the equations
of motion of the rings - desingularized to eliminate infinite self-induced
velocities - may be expressed as the 2-degree-of-freedom Hamiltonian system%
\begin{align}
\kappa _{k}\dot{s}_{k}& =4\kappa _{j}\kappa _{k}r_{j}r_{k}\left(
x_{k}-x_{j}\right) \int\nolimits_{0}^{\pi /2}\frac{\cos 2\sigma d\sigma }{%
\Delta _{12}^{3/2}},  \notag \\
&  \label{e1} \\
\kappa _{k}\dot{x}_{k}& =\left( \frac{\kappa _{k}^{2}}{2r_{k}}\right) \left[
\log \left( \frac{8r_{k}}{\delta }\right) -\gamma \right] +2\kappa
_{j}\kappa _{k}r_{j}\int\nolimits_{0}^{\pi /2}\frac{\left( r_{j}-r_{k}\cos
2\sigma \right) d\sigma }{\Delta _{12}^{3/2}},  \notag
\end{align}%
for $j,k=1,2$, with $j\neq k$, $\kappa _{k}:=\Gamma _{k}/2\pi $,
\begin{equation}
\Delta _{12}:=\left( r_{1}-r_{2}\right) ^{2}+\left( x_{1}-x_{2}\right)
^{2}+4r_{1}r_{2}\sin ^{2}\sigma ,  \label{e2}
\end{equation}%
$0<\delta \ll 1$ is a very small positive number representing the common
radius of the two rings in the desingularization procedure, and%
\begin{equation}
\gamma :=\frac{1}{2}\left( 1+\log 2+\int\nolimits_{0}^{\infty }e^{-\xi }\log
\xi d\xi \right) \cong 0.558,  \label{e3}
\end{equation}%
which is an artifact of assuming a Gaussian vorticity distribution in the
cores (\emph{i.e}. in the tubes of radius $\delta $) of the rings in the
desingularization. It should be remarked here that different
desingularization approaches produce slightly different parameters and
dynamical equations, but the overall qualitative aspects of the motion are
unaffected by such choices, and the differences in the quantitative
properties are quite small.

By introducing the Poisson bracket%
\begin{equation}
\left\{ f,g\right\} :=\sum\limits_{k=1}^{2}\frac{1}{\kappa_{k}}\left( \frac{%
\partial f}{\partial x_{k}}\frac{\partial g}{\partial s_{k}}-\frac{\partial f%
}{\partial s_{k}}\frac{\partial g}{\partial x_{k}}\right) ,  \label{e4}
\end{equation}
which is quite often employed in vortex dynamics, (1) can be recast in the
manifestly Hamiltonian form%
\begin{equation}
\dot{s}_{k}=\kappa_{k}^{-1}\partial_{x_{k}}H_{0}=\left\{ H_{0},s_{k}\right\}
,\;\dot{x}_{k}=-\kappa_{k}^{-1}\partial_{s_{k}}H_{0}=\left\{
H_{0},x_{k}\right\} ,\quad(1\leq k\leq2)  \label{e5}
\end{equation}
where the Hamiltonian function for the system is
\begin{align}
H_{0} & :=-\left\{ \sum\limits_{k=1}^{2}\kappa_{k}^{2}r_{k}\left[ \log\left(
\frac{8r_{k}}{\delta}\right) -\left( 1+\gamma\right) \right]
+4\kappa_{1}\kappa_{2}r_{1}r_{2}\int\nolimits_{0}^{\pi/2}\frac{\cos2\sigma
d\sigma}{\Delta_{12}^{1/2}}\right\}  \notag \\
& =-\left\{ \sum\limits_{k=1}^{2}\kappa_{k}^{2}r_{k}\left[ \log\left( \frac{%
8r_{k}}{\delta}\right) -\left( 1+\gamma\right) \right] +2\kappa
_{1}\kappa_{2}\left( r_{12+}+r_{12-}\right) \left[ F\left( \lambda
_{12}\right) -E\left( \lambda_{12}\right) \right] \right\} ,  \label{e6}
\end{align}
where%
\begin{equation}
r_{12\pm}^{2}:=\left( r_{1}\pm r_{2}\right) ^{2}+\left( z_{1}-z_{2}\right)
^{2},\;\lambda_{12}:=\frac{r_{12+}-r_{12-}}{r_{12+}+r_{12-}},  \label{e7}
\end{equation}
and $F$ and $G$ are, respectively, the complete elliptic integral of the
first and second kind given as
\begin{equation}
F\left( \lambda\right) :=\int\nolimits_{0}^{\pi/2}\frac{d\sigma}{\sqrt{%
1-\lambda^{2}\sin^{2}\sigma}}\;\mathrm{and\;}E\left( \lambda\right)
:=\int\nolimits_{0}^{\pi/2}\sqrt{1-\lambda^{2}\sin^{2}\sigma}d\sigma.
\label{e8}
\end{equation}
This use of elliptic integrals in the equations of motion is classical and
is often used (\emph{cf}. \cite{blkn1, HSC, lamb}). We note here that
(5)-(6) is (real) analytic on the following subset of $\mathbb{R}^{4}$,
which serves as the phase space of the system:%
\begin{equation*}
\boldsymbol{X}_{D}:=\left\{ \left( s_{1},s_{2},x_{1},x_{2}\right) \in\mathbb{%
R}^{4}:s_{1},s_{2}\geq0,\,\left( s_{1},x_{1}\right) \neq\left(
s_{2},x_{2}\right) \right\} .
\end{equation*}

The autonomous Hamiltonian system (5)-(6) has the following independent
integrals in involution (see \emph{e.g}. \cite{arn3, bb, blkn1, dbokVB,
mos1, new}):%
\begin{equation}
H_{0},\;G:=\sum\limits_{k=1}^{2}\kappa_{k}s_{k}=\sum\limits_{k=1}^{2}%
\kappa_{k}r_{k}^{2},  \label{e9}
\end{equation}
so it is completely integrable in the sense of Liouville-Arnold (\emph{%
LA-integrable}). We note here that a system of two or more coaxial rings -
even if one ring has a vortex strength of zero - does not have any
additional independent constants of motion in involution, as proved by
Bagrets \& Bagrets \cite{bb}. Consequently, the dynamics of three or more
rings is apt to include chaotic regimes in most cases.

\subsection{Full model including swirl}

We now immerse the two ring model in an ambient swirling flow. To be more
precise, we shall assume that the rings start and remain in a swirling flow
about the $x$-axis. In order to insure that the dynamical equations for the
full model including the coaxial vortex ring pair and the ambient flow still
have a Hamiltonian formulation, we shall choose our models for the ambient
swirling flow to be axisymmetric potential flows. A convenient way to
represent such swirling flows is via the use of cylindrical coordinates
having the $x$-axis (rather than the usual $z$-axis) as the axis of
symmetry. These coordinates naturally are defined as%
\begin{equation}
x=x,\;y=r\cos\theta,\;z=r\sin\theta,  \label{e10}
\end{equation}
where the angle $\theta$ is measured counterclockwise in the $yz$-plane
starting at zero along the positive $y$-axis. The general expression for the
ambient velocity in terms of these coordinates is%
\begin{equation}
\mathbf{v}_{a}=\left( \dot{x}_{a},\dot{y}_{a},\dot{z}_{a}\right) =\left(
u_{a},v_{a},w_{a}\right) =\left( \dot{x}_{a},\dot{r}_{a}\cos\theta-z_{a}\dot{%
\theta},\dot{r}_{a}\sin\theta+y_{a}\dot{\theta}\right) .  \label{e11}
\end{equation}
We assume for the ambient swirling flow that $\dot{x}_{a}$ is a quadratic
function of $s=r^{2}$ of the form
\begin{equation}
u_{a}\left( s\right) =\dot{x}_{a}\left( s\right) =-\alpha\left(
1+a_{1}s+a_{2}s^{2}\right) ,  \label{e12}
\end{equation}
where $\alpha>0$ and $a_{1}$ and $a_{2}$ are real constants to be chosen in
the sequel. To preserve axisymmetry, we assume that $\dot{\theta}$ is a
polynomial function of $s=r^{2}$ of degree two. Whence, we see that the
ambient velocity assumes the form%
\begin{equation}
\mathbf{v}_{a}=\left( \dot{x}_{a},\dot{y}_{a},\dot{z}_{a}\right) =\left(
\dot{x}_{a},\dot{r}_{a}\cos\theta-z_{a}\dot{\theta},\dot{r}_{a}\sin
\theta+y_{a}\dot{\theta}\right) =\left( -\alpha\left(
1+a_{1}s+a_{2}s^{2}\right) ,-z_{a}\dot{\theta}(s),y_{a}\dot{\theta}%
(s)\right) ,  \label{e13}
\end{equation}
where we shall only consider the following polynomial forms for the angular
velocity%
\begin{equation}
\dot{\theta}\left( s\right) =\Omega\left( 1+b_{1}s+b_{2}s^{2}\right) ,
\label{e14}
\end{equation}
where there is no loss of generality in assuming that $\Omega>0$, and $b_{1}$
and $b_{2}$ are real constants to be specified in the course of our
investigation.

Now we are in a position to formulate the full coaxial vortex ring pair +
swirling flow dynamics. Once again, the axisymmetry enables us to specify
the motion, modulo rotation about the $x$-axis, of the rings in terms of
their intersection points with the meridian half-plane $\mathfrak{H}$.
Observe that the \textquotedblleft modulo a rotation about the $x$%
-axis\textquotedblright\ is not really a problem in describing the evolution
of points on a ring. To see this, suppose we consider a point on a ring, say
the first ring, that is initially at $\left(
x_{1}(0),s_{1}(0),\theta_{1}(0)\right) $, and we know $x_{1}(t)$ and $%
s_{1}(t)$. Then we can find $\theta_{1}(t)$ from (14) and a simple
integration as%
\begin{equation}
\theta_{1}\left( t\right) :=\theta\left( s_{1}\left( t\right) \right)
=\theta_{1}(0)+\Omega\int\nolimits_{0}^{t}\left( 1+c_{1}s_{1}\left(
\tau\right) +c_{2}s_{1}^{2}\left( \tau\right) \right) d\tau,  \label{e15}
\end{equation}
and there is an analogous formula for $\theta_{2}(t):=\theta(s_{2}(t))$.

Adding the ambient swirl to (1) only changes the equations of motion by a
the addition of a simple polynomial function along the $x$-axes; namely, the
dynamical equations become

\begin{align}
\kappa _{k}\dot{s}_{k}& =4\kappa _{j}\kappa _{k}r_{j}r_{k}\left(
x_{k}-x_{j}\right) \int\nolimits_{0}^{\pi /2}\frac{\cos 2\sigma d\sigma }{%
\Delta _{12}^{3/2}},  \notag \\
&  \notag \\
\kappa _{k}\dot{x}_{k}& =-\kappa _{k}\alpha \left(
1+a_{1}s_{k}+a_{2}s_{k}^{2}\right) +\left( \frac{\kappa _{k}^{2}}{2r_{k}}%
\right) \left[ \log \left( \frac{8r_{k}}{\delta }\right) -\gamma \right]
+2\kappa _{j}\kappa _{k}r_{j}\int\nolimits_{0}^{\pi /2}\frac{\left(
r_{j}-r_{k}\cos 2\sigma \right) d\sigma }{\Delta _{12}^{3/2}},  \label{e16}
\end{align}%
having the following Hamiltonian analog of (5)-(6)

\begin{equation}
\dot{s}_{k}=\kappa_{k}^{-1}\partial_{x_{k}}H=\left\{ H,s_{k}\right\} ,\;\dot{%
x}_{k}=-\kappa_{k}^{-1}\partial_{s_{k}}H=\left\{ H,x_{k}\right\}
,\quad(1\leq k\leq2)  \label{e17}
\end{equation}
with Hamiltonian function
\begin{align}
H & :=\alpha\sum\limits_{k=1}^{2}\kappa_{k}\left[ s_{k}+\frac{a_{1}}{2}%
s_{k}^{2}+\frac{a_{2}}{3}s_{k}^{3}\right] +H_{0}=\alpha\sum\limits_{k=1}^{2}%
\kappa_{k}\left[ s_{k}+\frac{a_{1}}{2}s_{k}^{2}+\frac{a_{2}}{3}s_{k}^{3}%
\right] -  \notag \\
& \qquad\left\{ \sum\limits_{k=1}^{2}\kappa_{k}^{2}r_{k}\left[ \log\left(
\frac{8r_{k}}{\delta}\right) -\left( 1+\gamma\right) \right] +4\kappa
_{1}\kappa_{2}r_{1}r_{2}\int\nolimits_{0}^{\pi/2}\frac{\cos2\sigma d\sigma }{%
\Delta_{12}^{1/2}}\right\} =\alpha\sum\limits_{k=1}^{2}\kappa_{k}\left[
s_{k}+\frac{a_{1}}{2}s_{k}^{2}+\frac{a_{2}}{3}s_{k}^{3}\right] -  \notag \\
& \quad\quad\quad\quad\left\{ \sum\limits_{k=1}^{2}\kappa_{k}^{2}r_{k}\left[
\log\left( \frac{8r_{k}}{\delta}\right) -\left( 1+\gamma \right) \right]
+2\kappa_{1}\kappa_{2}\left( r_{12+}+r_{12-}\right) \left[ F\left(
\lambda_{12}\right) -E\left( \lambda_{12}\right) \right] \right\} .
\label{e18}
\end{align}
Observe that (17)-(18) also is LA-integrable inasmuch as it has the
following pair of independent invariants in involution:

\begin{equation}
H,\;G:=\sum\limits_{k=1}^{2}\kappa_{k}s_{k}=\sum\limits_{k=1}^{2}\kappa
_{k}r_{k}^{2}.  \label{e19}
\end{equation}

Our primary concern in the sequel is the affect of the ring dynamics - as
described by (17)-(18) - on the motion of passive fluid particles, which is
often described as the fluid kinematics.

\section{Dynamics of Model Kinematics}

The dynamics of passive fluid particles induced by the motion of the coaxial
vortex ring pair if the swirling flow characterized by\ (17)-(18) can be
described in two primary ways: Either indirectly and autonomously in terms
of a restricted three ring problem in which we treat a passive fluid
particle as belonging to a third ring of zero strength, or directly as a
time-dependent system of ordinary differential equations in which we first
solve for the dynamics of (16) in order to describe its affect on a passive
fluid particle. Each of these approaches has its advantages, and we shall
find it convenient to develop both of them in what follows.

\subsection{Passive third ring approach to the kinematics}

In this approach, we actually consider a three ring problem with two of the
rings as above with their given nonzero strengths, $\Gamma_{1}$ and $%
\Gamma_{2}$, and a third (advected) ring of strength $\Gamma=0$. Once again,
the axisymmetry implies that the motion of the rings is completely
determined by that of their respective intersection points, $(x_{1}.y_{1})$,
$(x_{2},y_{2}) $ and $(x,y)$ with the half-plane $\mathfrak{H}$. The
equations of motion of these three points are analogous to (16), and can be
written as(\emph{cf.} \cite{bb, blkn1, blkn2, BWC, BVS, BB})

\begin{align}
\kappa _{k}\dot{s}_{k}& =4\kappa _{j}\kappa _{k}r_{j}r_{k}\left(
x_{k}-x_{j}\right) \int\nolimits_{0}^{\pi /2}\frac{\cos 2\sigma d\sigma }{%
\Delta _{12}^{3/2}},  \notag \\
\kappa _{k}\dot{x}_{k}& =-\kappa _{k}\alpha \left(
1+a_{1}s_{k}+a_{2}s_{k}^{2}\right) +\left( \frac{\kappa _{k}^{2}}{2r_{k}}%
\right) \left[ \log \left( \frac{8r_{k}}{\delta }\right) -\gamma \right]
+2\kappa _{j}\kappa _{k}r_{j}\int\nolimits_{0}^{\pi /2}\frac{\left(
r_{j}-r_{k}\cos 2\sigma \right) d\sigma }{\Delta _{12}^{3/2}},  \notag \\
\dot{s}& =r\sum\limits_{k=1}^{2}\kappa _{k}r_{k}\left( x-x_{k}\right)
\int\nolimits_{0}^{\pi /2}\frac{\cos 2\sigma d\sigma }{\Delta _{k}^{3/2}},
\notag \\
\dot{x}& =-\alpha \left( 1+a_{1}s+a_{2}s^{2}\right)
+2\sum\limits_{k=1}^{2}\kappa _{k}r_{k}\int\nolimits_{0}^{\pi /2}\frac{%
\left( r_{k}-r\cos 2\sigma \right) d\sigma }{\Delta _{k}^{3/2}},  \label{e20}
\end{align}%
where

\begin{equation}
\Delta_{k}:=\left( r-r_{k}\right) ^{2}+\left( x-x_{k}\right)
^{2}+4rr_{k}\sin^{2}\sigma.\quad(1\leq k\leq2)  \label{e21}
\end{equation}

This system can be recast in what amounts to a kind of \emph{piecewise
Hamiltonian} form. In particular, we may rewrite (20) as

\begin{align}
\dot{s}_{k} & =\kappa_{k}^{-1}\partial_{x_{k}}H=\left\{ H,s_{k}\right\} ,\;%
\dot{x}_{k}=-\kappa_{k}^{-1}\partial_{s_{k}}H_{\ast}=\left\{ H,x_{k}\right\}
,\quad(1\leq k\leq2)  \notag \\
\dot{s} & =\partial_{x}H_{\ast}=\left\{ H_{\ast},s\right\} _{\ast},\;\dot{x}%
=-\partial_{s}H_{\ast}=\left\{ H_{\ast},x\right\} _{\ast},  \label{e22}
\end{align}
where $\{\cdot,\cdot\}$ is as in (4), $H$ is as in (18), and

\begin{align}
H_{\ast} & :=H+\alpha\left[ s+\frac{a_{1}}{2}s^{2}+\frac{a_{2}}{3}s^{3}%
\right] -4r\sum\limits_{k=1}^{2}\kappa_{k}r_{k}\int\nolimits_{0}^{\pi /2}%
\frac{\cos2\sigma d\sigma}{\Delta_{k}^{1/2}}=  \notag \\
& \alpha\left\{ \left[ s+\frac{a_{1}}{2}s^{2}+\frac{a_{2}}{3}s^{3}\right]
+\sum\limits_{k=1}^{2}\kappa_{k}\left[ s_{k}+\frac{a_{1}}{2}s_{k}^{2}+\frac{%
a_{2}}{3}s_{k}^{3}\right] \right\} -\left\{
\sum\limits_{k=1}^{2}\kappa_{k}^{2}r_{k}\left[ \log\left( \frac{8r_{k}}{%
\delta}\right) -\left( 1+\gamma\right) \right] +\right.  \notag \\
& \qquad\qquad\left. 4\kappa_{1}\kappa_{2}r_{1}r_{2}\int\nolimits_{0}^{\pi/2}%
\frac{\cos2\sigma d\sigma}{\Delta_{12}^{1/2}}+4r\sum\limits_{k=1}^{2}%
\kappa_{k}r_{k}\int\nolimits_{0}^{\pi/2}\frac{\cos2\sigma d\sigma}{%
\Delta_{k}^{1/2}}\right\}  \notag \\
& =\alpha\left\{ \left[ s+\frac{a_{1}}{2}s^{2}+\frac{a_{2}}{3}s^{3}\right]
+\sum\limits_{k=1}^{2}\kappa_{k}\left[ s_{k}+\frac{a_{1}}{2}s_{k}^{2}+\frac{%
a_{2}}{3}s_{k}^{3}\right] \right\} -\left\{
\sum\limits_{k=1}^{2}\kappa_{k}^{2}r_{k}\left[ \log\left( \frac{8r_{k}}{%
\delta}\right) -\left( 1+\gamma\right) \right] +\right.  \notag \\
& \qquad\qquad\quad\left. 2\kappa_{1}\kappa_{2}\left( r_{12+}+r_{12-}\right)
\left[ F\left( \lambda_{12}\right) -E\left( \lambda _{12}\right) \right]
+2\sum\limits_{k=1}^{2}\kappa_{k}\left( r_{k+}+r_{k-}\right) \left[ F\left(
\lambda_{k}\right) -E\left( \lambda _{k}\right) \right] \right\} ,
\label{e23}
\end{align}
where%
\begin{equation}
r_{k\pm}^{2}:=\left( r_{k}\pm r\right) ^{2}+\left( z_{k}-z\right)
^{2},\;\lambda_{k}:=\frac{r_{k+}-r_{k-}}{r_{k+}+r_{k-}},\quad(1\leq k\leq2)
\label{e24}
\end{equation}
$F$ and $E$ are as in (8), and the new Poisson bracket is actually the
standard one defined as

\begin{equation}
\left\{ f,g\right\} _{\ast}:=\left( \frac{\partial f}{\partial x}\frac{%
\partial g}{\partial s}-\frac{\partial f}{\partial s}\frac{\partial g}{%
\partial x}\right) .  \label{e25}
\end{equation}
The system (22)-(23) is analytic on the six-dimensional symplectic space%
\begin{equation*}
\boldsymbol{X}_{DK}:=\left\{ \left( s_{1},s_{2},x_{1},x_{2},s,x\right) \in%
\mathbb{R}^{6}:s_{1},s_{2},s\geq0,\left( s_{1},s_{2},x_{1},x_{2}\right) \in%
\boldsymbol{X}_{D},\,\left( s,x\right) \neq\left( s_{k},x_{k}\right)
,\,1\leq k\leq2\,\right\} ,
\end{equation*}
and clearly shares the following two independent constants of motion with
(16):%
\begin{equation}
H,\;G:=\sum\limits_{k=1}^{2}\kappa_{k}s_{k}=\sum\limits_{k=1}^{2}\kappa
_{k}r_{k}^{2},  \label{e26}
\end{equation}
but no additional integrals in involution as demonstrated in \cite{bb}. From
this we infer that the motion of the advected ring can become chaotic, which
is consistent with the apparent chaotic motion found in vortex bubbles \cite%
{BVS, BREU, BB, GBS, LP, serbon, sotir1, sotir2}.

\subsection{Direct formulation of kinematics}

In the direct method, we obtain a time-dependent system of differential
equations for the motion of a general passive fluid particle by assuming
that the dynamics of the pair of coaxial rings has already been determined;
\emph{i.e}., we have the solution $(s_{1}(t),x_{1}(t)),(s_{2}(t),x_{2}(t))$
of (17)-(18), which means that we also know $\theta_{1}(t):=\theta(s_{1}(t))$
and $\theta_{2}(t):=\theta(s_{2}(t))$ owing to (15).

We could start with the Biot-Savart law as was done in Blackmore \& Knio
\cite{blkn1} and derive the system of differential equations governing the
motion of points $\mathbf{x}=(x,y,z)\in \mathbb{R}^{3}$. However, we
simplify things by taking advantage of the fact that the rotational motion
of the particle about the $x$-axis is completely determined via (14)-(16) by
the motion of the intersection point with $\mathfrak{H}$ of the advected
circular ring containing the particle. Not only does this reduce the
dimension of the (continuous) dynamical system from three to two, but it
also has the advantage of having a Hamiltonian formulation, albeit a
time-dependent one. Taking this approach, it is easy to see that the desired
system is already at hand - namely from the last two equations of (20) with
the positions of the two rings specified as functions of time. In
particular, the equations of motion of the intersection points of the
advected ring with $\mathfrak{H}$ are%
\begin{align}
\dot{s}& =4r\sum\limits_{k=1}^{2}\kappa _{k}r_{k}\left( x-x_{k}\left(
t\right) \right) \int\nolimits_{0}^{\pi /2}\frac{\cos 2\sigma d\sigma }{%
\Delta _{k}^{3/2}\left( t\right) },  \notag \\
\dot{x}& =-\alpha \left( 1+a_{1}s+a_{2}s^{2}\right)
+2\sum\limits_{k=1}^{2}\kappa _{k}r_{k}(t)\int\nolimits_{0}^{\pi /2}\frac{%
\left( r_{k}\left( t\right) -r\cos 2\sigma \right) d\sigma }{\Delta
_{k}^{3/2}\left( t\right) },  \label{e27}
\end{align}%
where,
\begin{equation}
\Delta _{k}\left( t\right) :=\left( r-r_{k}\left( t\right) \right)
^{2}+\left( x-x_{k}\left( t\right) \right) ^{2}+4rr_{k}\left( t\right) \sin
^{2}\sigma .\quad (1\leq k\leq 2)  \label{e28}
\end{equation}

The system (27) is (time-dependent) Hamiltonian since it can be expressed in
the form%
\begin{equation}
\dot{s}=\partial_{x}\mathcal{H}\left( s,x,t\right) ,\quad\dot{x}%
=-\partial_{s}\mathcal{H}\left( s,x,t\right) ,  \label{e29}
\end{equation}
where the Hamiltonian function is%
\begin{align}
\mathcal{H}\left( s,x,t\right) & :=\alpha\left( s+\frac{a_{1}}{2}s^{2}+\frac{%
a_{3}}{3}s^{3}\right) -4r\sum\limits_{k=1}^{2}\kappa_{k}r_{k}\left( t\right)
\int\nolimits_{0}^{\pi/2}\frac{\cos2\sigma d\sigma }{\Delta_{k}^{1/2}\left(
t\right) }  \notag \\
& =\alpha\left( s+\frac{a_{1}}{2}s^{2}+\frac{a_{3}}{3}s^{3}\right)
-2\sum\limits_{k=1}^{2}\kappa_{k}\left( r_{k+}\left( t\right) +r_{k-}\left(
t\right) \right) \left[ F\left( \lambda_{k}(t)\right) -E\left(
\lambda_{k}(t)\right) \right] ,  \label{e30}
\end{align}
and we have used the notation in (24) interpreted as (known) functions of
time. We note that this system is analytic on the phase space%
\begin{equation*}
\boldsymbol{X}_{K}:=\left\{ \left( s,x,t\right) \in\mathbb{R}%
^{3}:s\geq0,\,\left( s,x\right) \neq\left( s_{k}(t),x_{k}(t)\right) ,\,1\leq
k\leq2\right\} .
\end{equation*}
Observe that it is clear from either (20) or (27) that if the passive
particle starts on the $x$-axis ($s=0$), it remains there for all time. To
be more precise, the hypersurface $s=0$ is an invariant hypersurface of the
six-dimensional phase space for (20), and it is an invariant line in the
trajectory space associated to (27). From a fluid dynamics perspective, the
points on $s=0$ for which $\dot{x}=0$ are of particular importance since
they correspond to stagnation points on the axis of symmetry, which one
would expect to define the leading and trailing edge (points) of a B-type
vortex breakdown structure. We shall have much more to say about this in the
course of our analysis, which follows. Another observation that will prove
useful for the analysis, is that the half-plane plays a very natural role as
a transversal for a fairly obvious Poincar\'{e} map - a map that conveys a
great deal of useful dynamical information. As a preview, we give a brief
description of this map. Let $(s(t),x(t))$ be the unique solution of (27)
starting at $(s_{0},x_{0})\in\mathfrak{H}$, with $s_{0}>0$. Then the
trajectory in $\mathbb{R}^{3}$ winds around the $x$-axis until at a first $%
t_{i}>0$, it intersects $\mathfrak{H}$ again. This establishes a smooth
mapping of $\mathfrak{H}$ into itself via the formula $\mathit{\Pi}\left(
(s_{0},x_{0})\right) :=(s(t_{i}),x(t_{i}))$. This handy mapping shall be
investigated further in the sequel.

\section{Analysis of the Ring Dynamics}

As for the ring dynamics, governed by (17)-(18), it was pointed out in the
preceding section that it is LA-integrable, a property often referred to as
integrable by quadratures. Therefore, at least in a theoretical sense, we
can integrate the dynamical equations (17)-(18). It appears, however, that
this integration cannot generally be carried out in a practical way in order
to find a tractable closed form solution for the equations of motion, except
perhaps in some very special cases.

\subsection{Stationary points}

We begin with an investigation of the fixed points of the system, which we
consider in the form of (16). Solutions of $\dot{s}_{1}=\dot{s}_{2}=0$ are
easy to characterize. This follows from the readily verifiable fact that the
integral
\begin{equation}
I:=\int\nolimits_{0}^{\pi/2}\frac{\cos2\sigma d\sigma}{\Delta_{12}^{3/2}}
\label{e31}
\end{equation}
is positive for all choices of $r_{1},r_{2}>0,x_{1}$ and $x_{2}$ ( except in
a limiting case where one ring collapses to a point with its core radius
decreasing in a special way). Thus a necessary and sufficient condition for $%
\dot{s}_{1}=\dot{s}_{2}=0$ is that $x_{1}=x_{2}$. This common value shall be
denoted as $\xi:=x_{1}=x_{2}$.

Naturally, it remains to find solutions of $\ \dot{x}_{1}=\dot{x}_{2}=0$,
and to this end we introduce some convenient simplifications. There is no
loss of generality in assuming
\begin{equation}
1=\kappa _{1}\leq \kappa :=\kappa _{2}\;\mathrm{and}\;0<r_{1}<r_{2}
\label{e32}
\end{equation}%
since this can be attained by a rescaling of the time variable $t$ and
reordering the variables, if necessary. In order to simplify matters, we
shall assume in (12) that $a_{1}=a_{2}=0$, so the equations we need to solve
take the form%
\begin{align}
0& =\dot{x}_{1}=-\alpha +\frac{1}{2r_{1}}\left[ \log \left( \frac{8r_{1}}{%
\delta }\right) -\gamma \right] +2\kappa r_{2}\int\nolimits_{0}^{\pi /2}%
\frac{\left( r_{2}-r_{1}\cos 2\sigma \right) d\sigma }{\hat{\Delta}%
_{12}^{3/2}},  \notag \\
0& =\dot{x}_{2}=-\alpha +\frac{\kappa }{2r_{2}}\left[ \log \left( \frac{%
8r_{2}}{\delta }\right) -\gamma \right] +2r_{1}\int\nolimits_{0}^{\pi /2}%
\frac{\left( r_{1}-r_{2}\cos 2\sigma \right) d\sigma }{\hat{\Delta}%
_{12}^{3/2}},  \label{e33}
\end{align}%
where
\begin{equation}
\mathring{\Delta}_{12}:=\left( r_{1}-r_{2}\right) ^{2}+4r_{1}r_{2}\sin
^{2}\sigma .  \label{e34}
\end{equation}%
It is convenient to display the full set of equations defining the fixed
points, which owing to (16) and our assumptions can be written as%
\begin{align}
\dot{s}_{1}& =\Phi _{1}\left( r_{1},r_{2},x_{1},x_{2};\alpha ,\kappa \right)
:=4\kappa r_{1}r_{2}\left( x_{1}-x_{2}\right) \int\nolimits_{0}^{\pi /2}%
\frac{\cos 2\sigma d\sigma }{\Delta _{12}^{3/2}}=0,  \notag \\
\dot{s}_{2}& =\Phi _{2}\left( r_{1},r_{2},x_{1},x_{2};\alpha ,\kappa \right)
:=4r_{1}r_{2}\left( x_{2}-x_{1}\right) \int\nolimits_{0}^{\pi /2}\frac{\cos
2\sigma d\sigma }{\Delta _{12}^{3/2}}=0,  \notag \\
\dot{x}_{1}& =\Psi _{1}\left( r_{1},r_{2},x_{1},x_{2};\alpha ,\kappa \right)
:=-\alpha +\frac{1}{2r_{1}}\left[ \log \left( \chi r_{1}\right) -\gamma %
\right] +2\kappa r_{2}\int\nolimits_{0}^{\pi /2}\frac{\left( r_{2}-r_{1}\cos
2\sigma \right) d\sigma }{\Delta _{12}^{3/2}}=0,  \notag \\
\dot{x}_{2}& =\Psi _{2}\left( r_{1},r_{2},x_{1},x_{2};\alpha ,\kappa \right)
:=-\alpha +\frac{\kappa }{2r_{2}}\left[ \log \left( \chi r_{2}\right)
-\gamma \right] +2r_{1}\int\nolimits_{0}^{\pi /2}\frac{\left(
r_{1}-r_{2}\cos 2\sigma \right) d\sigma }{\Delta _{12}^{3/2}}=0,  \label{e35}
\end{align}%
where $\chi \gg 1$ is a multiple of the curvature of circular cross-sections
of the boundary of the (virtual) core of each coaxial vortex ring, and is
defined as%
\begin{equation}
\chi :=\frac{8}{\delta }.  \label{e36}
\end{equation}%
We note that, strictly speaking, the functions in (35) also depend on the
parameter $\chi $; however, given the \textquotedblleft virtual
nature\textquotedblright\ of the definition of $\chi $, we choose to fix its
value at 1,000 and suppress it as a coordinate of parameter space in the
sequel. As the first two equations of this system require that $%
x_{1}=x_{2}:=\xi $, this solution reduces to the following pair of \ $\xi $%
-independent \emph{radii equations}:%
\begin{align}
\hat{\Psi}_{1}\left( r_{1},r_{2};\alpha ,\kappa \right) & :=-\alpha +\frac{1%
}{2r_{1}}\left[ \log \left( \chi r_{1}\right) -\gamma \right] +2\kappa
r_{2}\int\nolimits_{0}^{\pi /2}\frac{\left( r_{2}-r_{1}\cos 2\sigma \right)
d\sigma }{\mathring{\Delta}_{12}^{3/2}}=0,  \notag \\
\hat{\Psi}_{2}\left( r_{1},r_{2};\alpha ,\kappa \right) & :=-\alpha +\frac{%
\kappa }{2r_{2}}\left[ \log \left( \chi r_{2}\right) -\gamma \right]
+2r_{1}\int\nolimits_{0}^{\pi /2}\frac{\left( r_{1}-r_{2}\cos 2\sigma
\right) d\sigma }{\mathring{\Delta}_{12}^{3/2}}=0.  \label{e37}
\end{align}

\subsubsection{Solution of radii equations}

It can be shown via a straightforward argument, in which the function%
\begin{equation}
\psi\left( r\right) :=\left( 2r\right) ^{-1}\left[ \log\left( \chi r\right)
-\gamma\right]  \label{e38}
\end{equation}
plays a key role, that the radii equations (37) have precisely four
solutions $(r_{1},r_{2})$ when $\alpha$ lies between the value one and a not
too large upper bound. In these four solutions, there is one with both
coordinates very small, one with both coordinates significantly larger than
in the first solution, and two solutions with one coordinate small and the
other quite a bit larger. We leave the proof of this to the reader, and only
give a sketch of the reasoning involved. Toward this end, it is useful to
note some basic properties of the graph of the function $\psi$ restricted
naturally to $r>0$. It is easy to see that $\psi$ has a unique zero at $%
r_{0}:=e^{\gamma}\chi ^{-1}$ and its derivative $\psi^{\prime}$ also has a
unique zero at $r=\rho:=e^{\left( \gamma+1\right) }\chi^{-1}$ at which the
function attains its maximum of $\psi\left( \rho\right) =2^{-1}e^{-\left(
\gamma+1\right) }\chi=\rho^{-1}$. Moreover, $\psi^{\prime}\left( r\right) $
is positive for $0<r<\rho$ and negative for $\rho<r<\infty$, and $\psi\left(
r\right) \rightarrow0$ as $r\rightarrow\infty$.

For our sketch of the solution of the radii equations, we shall find it
instructive to ignore the integral terms in (37) thereby obtaining the \emph{%
reduced radii equations}

\begin{align}
\breve{\Psi}_{1}\left( r_{1},r_{2};\alpha,\kappa\right) & :=-\alpha +\frac{1%
}{2r_{1}}\left[ \log\left( \chi r_{1}\right) -\gamma\right] =0,  \notag \\
\breve{\Psi}_{2}\left( r_{1},r_{2};\alpha,\kappa\right) & :=-\alpha +\frac{%
\kappa}{2r_{2}}\left[ \log\left( \chi r_{2}\right) -\gamma\right] =0.
\label{e39}
\end{align}
It turns out that the solutions of the full equations (37) are completely
analogous to those of (39), and can be obtained rather directly using
identical methods. Accordingly it follows from the form of (39) that we
ought to restrict $\alpha$ in both (39) and (37) so that $0<\alpha<\rho^{-1}$%
, which is not much of a restriction at all since $\rho^{-1}$ is so large.
It is convenient to recast (39) in the form
\begin{align}
-\left( 2\alpha r_{1}+\gamma\right) +\log\left( \chi r_{1}\right) & =0,
\notag \\
-\left( 2\alpha r_{2}+\kappa\gamma\right) +\kappa\log\left( \chi
r_{2}\right) & =0.  \label{e40}
\end{align}
It is easy to show that for $1<\alpha<\rho^{-1}$, each of (40) has a small
(positive) solution and a larger solution to the left and right,
respectively, of $\rho$. Let us denote these solutions, which are functions
of $(\alpha,\kappa)$, as $\tilde{r}_{1_{a}},\tilde{r}_{1_{b}},\tilde{r}%
_{2_{a}}$ and $\tilde{r}_{2_{b}}$, and note that is a simple matter to
verify that they satisfy the following properties:%
\begin{equation*}
r_{0}<\tilde{r}_{2_{a}}<\tilde{r}_{1_{a}}<\rho<\tilde{r}_{1_{b}}<\tilde {r}%
_{2_{b}},
\end{equation*}%
\begin{equation*}
\tilde{r}_{1_{a}}(\alpha,\kappa)-\tilde{r}_{2_{a}}(\alpha,\kappa),\tilde {r}%
_{2_{b}}(\alpha,\kappa)-\tilde{r}_{1_{b}}(\alpha,\kappa)\;\mathrm{and\;}%
\tilde{r}_{2_{b}}(\alpha,\kappa)-\tilde{r}_{1_{a}}(\alpha,\kappa)
\end{equation*}
are decreasing functions of $\alpha$ for every $\kappa>1$ and increasing
functions of $\kappa$ for each $\alpha>1$, and finally%
\begin{equation*}
\tilde{r}_{2_{a}}(\alpha,\kappa)\downarrow r_{0}\;\mathrm{as}\;\kappa
\uparrow\infty
\end{equation*}
for each $\alpha>1$. We also note that it is easy to show that $\tilde {r}%
_{1_{a}}$ can be obtained as the limit of the following iterative scheme in
which the first equation of (40) has been reformulated as a fixed point
problem%
\begin{equation}
r_{1}^{(n+1)}=\frac{1}{\chi}\exp\left( 2\alpha r_{1}^{(n)}+\gamma\right) ,
\label{e41}
\end{equation}
where the initial iterate $r_{1}^{(0)}$ can be chosen as any value between
zero and $\rho$. Moreover, $\tilde{r}_{1_{+}}$ can be obtained as the limit
of the following iterative scheme in which the first equation of (40) is
reformulated as the (inverse) fixed point problem
\begin{equation}
r_{1}^{(n+1)}=\frac{1}{2\alpha}\left[ \log\left( \chi r_{1}^{(n)}\right)
-\gamma\right] ,  \label{e42}
\end{equation}
where the initial iterate $r_{1}^{(0)}$ in this case can be chosen to be any
reasonably large number greater than $\rho$; for example, $r_{1}^{(0)}=10$
is always a suitable choice. Naturally, there are completely analogous
schemes for finding $\tilde{r}_{2_{a}}$ and $\tilde{r}_{2_{b}}$ using the
second of the equations (40).

These observations about the reduced form (39) of (37) provide the basic
tools for solving (37). More precisely, we find that although the additional
integrals in (37) break some of the symmetry inherent in (39), there are all
together four solutions, which we denote as $(\hat{r}_{1(\mathrm{I})},\hat{r}%
_{2(\mathrm{I})})=(\hat{r}_{1(\mathrm{I})}\left( \alpha ,\kappa \right) ,%
\hat{r}_{2(\mathrm{I})}\left( \alpha ,\kappa \right) )$, $(\hat{r}_{1(%
\mathrm{II})},\hat{r}_{2(\mathrm{II})})=(\hat{r}_{1(\mathrm{II})}\left(
\alpha ,\kappa \right) ,\hat{r}_{2(\mathrm{II})}\left( \alpha ,\kappa
\right) )$, $(\hat{r}_{1(\mathrm{III})},\hat{r}_{2(\mathrm{III})})=(\hat{r}%
_{1(\mathrm{III})}\left( \alpha ,\kappa \right) ,\hat{r}_{2(\mathrm{III}%
)}\left( \alpha ,\kappa \right) )$, and $(\hat{r}_{1(\mathrm{IV})},\hat{r}%
_{2(\mathrm{IV})})=(\hat{r}_{1(\mathrm{IV})}\left( \alpha ,\kappa \right) ,%
\hat{r}_{2(\mathrm{IV})}\left( \alpha ,\kappa \right) )$, and for which the
following properties - somewhat mirroring those of the solutions of (39) -
hold:

\begin{align}
0& <\hat{r}_{1(\mathrm{I})}<\hat{r}_{2(\mathrm{I})}\simeq 4\hat{r}_{1(%
\mathrm{I})};\;0<\hat{r}_{2(\mathrm{II})}\ll \hat{r}_{1(\mathrm{II})}\simeq
\left( 2\times 10^{3}\right) \hat{r}_{2(\mathrm{II})};  \notag \\
& 0<\hat{r}_{1(\mathrm{III})}\ll \hat{r}_{2(\mathrm{III})}\simeq \left(
5\times 10^{2}\right) \hat{r}_{1(\mathrm{III})};\;0<\left( 10^{-2}\right)
\hat{r}_{1(\mathrm{I})}\simeq \hat{r}_{1(\mathrm{IV})}<\hat{r}_{2(\mathrm{IV}%
)}\simeq 3\hat{r}_{1(\mathrm{IV})};  \label{e43} \\
& \mathrm{and\;}\hat{r}_{2(\mathrm{I})}-\hat{r}_{1(\mathrm{I})},\;\hat{r}_{1(%
\mathrm{II})}-\hat{r}_{2(\mathrm{II})},\;\hat{r}_{2(\mathrm{III})}-\hat{r}%
_{1(\mathrm{III})},\;\mathrm{and}\;\hat{r}_{2(\mathrm{IV})}-\hat{r}_{1(%
\mathrm{IV})}  \notag
\end{align}%
are all decreasing functions of $\alpha $ for every $\kappa >1$ and
increasing functions of $\kappa $ for each $\alpha >1$. We also note that
all four solutions can be obtained iteratively using the procedures
described in (41) and (42), although considerably more care must be
exercised with regard to selecting the initial point in the process. In
particular, the solution $(\hat{r}_{1(\mathrm{III})},\hat{r}_{2(\mathrm{III}%
)})$ may be obtained via Picard iteration in the form%
\begin{align}
r_{1}^{(n+1)}& =\frac{1}{\chi }\exp \left[ \left( 2\alpha r_{1}^{(n)}+\gamma
\right) -4\kappa r_{1}^{(n)}r_{2}^{(n)}\int\nolimits_{0}^{\pi /2}\frac{%
\left( r_{2}^{(n)}-r_{1}^{(n)}\cos 2\sigma \right) d\sigma }{\mathring{\Delta%
}_{12}^{3/2}\left( r_{1}^{(n)},r_{2}^{(n)}\right) }\right] ,  \notag \\
r_{2}^{(n+1)}& =\frac{1}{2\alpha }\left\{ \kappa \left[ \log \left( \chi
r_{2}^{(n)}\right) -\gamma \right] +4r_{1}^{(n)}r_{2}^{(n)}\int%
\nolimits_{0}^{\pi /2}\frac{\left( r_{1}^{(n)}-r_{2}^{(n)}\cos 2\sigma
\right) d\sigma }{\mathring{\Delta}_{12}^{3/2}\left(
r_{1}^{(n)},r_{2}^{(n)}\right) }\right\} ,  \label{e44}
\end{align}%
with $\left( r_{1}^{(0)},r_{2}^{(0)}\right) =\left( 0.001,1.0\right) $, and
then
\begin{equation*}
\left( r_{1}^{(n)},r_{2}^{(n)}\right) \rightarrow \left( \hat{r}_{1(\mathrm{%
III})},\hat{r}_{2(\mathrm{III})}\right) \;\mathrm{as}\;n\rightarrow \infty .
\end{equation*}

\subsection{Stability of stationary points}

We have now identified a whole 1-parameter family of stationary points for
(16) that can be represented as in terms of the parameter $\xi $ as
\begin{align}
\mathfrak{L}& =\mathfrak{L}\left( \xi \right) :=\left\{ q=\left(
s_{1},s_{2},x_{1},x_{2}\right) :s_{1}=\hat{s}_{1(\mathrm{I})}:=\hat{r}_{1(%
\mathrm{I})}^{2},s_{1}=\hat{s}_{1(\mathrm{II})}:=\hat{r}_{1(\mathrm{II}%
)}^{2},s_{1}=\hat{s}_{1(\mathrm{III})}:=\hat{r}_{1(\mathrm{III})}^{2},\right.
\notag \\
& \quad \quad \quad \quad \quad \;\;\left. s_{1}=\hat{s}_{1(\mathrm{IV})}:=%
\hat{r}_{1(\mathrm{IV})}^{2},s_{2}=\hat{s}_{2(\mathrm{I})}:=\hat{r}_{2(%
\mathrm{I})}^{2},s_{2}=\hat{s}_{2(\mathrm{II})}:=\hat{r}_{2(\mathrm{II}%
)}^{2},\,s_{2}=\hat{s}_{2(\mathrm{III})}:=\hat{r}_{2(\mathrm{III}%
)}^{2},\right.  \label{eq45} \\
& \qquad \qquad \qquad \left. s_{2}=\hat{s}_{2(\mathrm{IV})}:=\hat{r}_{2(%
\mathrm{IV})}^{2},x_{1}=x_{2}=\xi ,\;\xi \in \mathbb{R}\right\} .  \notag
\end{align}%
To perform a standard linear stability analysis for a typical $q\in
\mathfrak{L}$, we compute the derivative (matrix) for $\mathfrak{X}%
_{12}:=\left( \Phi _{1},\Phi _{2},\Psi _{1},\Psi _{2}\right) $; namely,%
\begin{equation}
\mathfrak{X}_{12}^{\prime }\left( q\right) :=\left(
\begin{array}{cccc}
\partial _{s_{1}}\Phi _{1}\left( q\right) & \partial _{s_{2}}\Phi _{1}\left(
q\right) & \partial _{x_{1}}\Phi _{1}\left( q\right) & \partial _{x_{2}}\Phi
_{1}\left( q\right) \\
\partial _{s_{1}}\Phi _{2}\left( q\right) & \partial _{s_{2}}\Phi _{2}\left(
q\right) & \partial _{x_{1}}\Phi _{2}\left( q\right) & \partial _{x_{2}}\Phi
_{2}\left( q\right) \\
\partial _{s_{1}}\Psi _{1}\left( q\right) & \partial _{s_{2}}\Psi _{1}\left(
q\right) & \partial _{x_{1}}\Psi _{1}\left( q\right) & \partial _{x_{2}}\Psi
_{1}\left( q\right) \\
\partial _{s_{1}}\Psi _{2}\left( q\right) & \partial _{s_{2}}\Psi _{2}\left(
q\right) & \partial _{x_{1}}\Psi _{2}\left( q\right) & \partial _{x_{2}}\Psi
_{2}\left( q\right)%
\end{array}%
\right) .  \label{e46}
\end{equation}%
In aid of this, we compute that%
\begin{align}
\partial _{s_{1}}\Phi _{1}\left( q\right) & =0,\;\partial _{s_{2}}\Phi
_{1}\left( q\right) =0,\;\partial _{x_{1}}\Phi _{1}\left( q\right) =4\kappa
\hat{r}_{1}\hat{r}_{2}\int\nolimits_{0}^{\pi /2}\frac{\cos 2\sigma d\sigma }{%
\mathring{\Delta}_{12}^{3/2}\left( \hat{r}_{1},\hat{r}_{2}\right) }%
,\;\partial _{x_{2}}\Phi _{1}\left( q\right) =-4\kappa \hat{r}_{1}\hat{r}%
_{2}\int\nolimits_{0}^{\pi /2}\frac{\cos 2\sigma d\sigma }{\mathring{\Delta}%
_{12}^{3/2}\left( \hat{r}_{1},\hat{r}_{2}\right) },  \notag \\
\partial _{s_{1}}\Phi _{2}\left( q\right) & =0,\;\partial _{s_{2}}\Phi
_{2}\left( q\right) =0,\;\partial _{x_{1}}\Phi _{2}\left( q\right) =-4\hat{r}%
_{1}\hat{r}_{2}\int\nolimits_{0}^{\pi /2}\frac{\cos 2\sigma d\sigma }{%
\mathring{\Delta}_{12}^{3/2}\left( \hat{r}_{1},\hat{r}_{2}\right) }%
,\;\partial _{x_{2}}\Phi _{2}\left( q\right) =4\hat{r}_{1}\hat{r}%
_{2}\int\nolimits_{0}^{\pi /2}\frac{\cos 2\sigma d\sigma }{\mathring{\Delta}%
_{12}^{3/2}\left( \hat{r}_{1},\hat{r}_{2}\right) },  \notag \\
\partial _{s_{1}}\Psi _{1}\left( q\right) & =-\frac{1}{2\hat{r}_{1}}\left\{
\frac{1}{2\hat{r}_{1}^{2}}\left[ \log \left( \chi \hat{r}_{1}\right) -\left(
1+\gamma \right) \right] -\right.  \notag \\
& \quad \left. 2\kappa \hat{r}_{2}\int\nolimits_{0}^{\pi /2}\frac{\left(
\hat{\Delta}_{12}\left( \hat{r}_{1},\hat{r}_{2}\right) \cos 2\sigma +3\left(
\hat{r}_{2}-\hat{r}_{1}\cos 2\sigma \right) \left[ \left( \hat{r}_{1}-\hat{r}%
_{2}\right) +2\hat{r}_{2}\sin ^{2}\sigma \right] \right) d\sigma }{\mathring{%
\Delta}_{12}^{5/2}\left( \hat{r}_{1},\hat{r}_{2}\right) }\right\} ,  \notag
\\
\partial _{s_{2}}\Psi _{1}\left( q\right) & =\frac{\kappa }{\hat{r}_{2}}%
\int\nolimits_{0}^{\pi /2}\frac{\left( \left( 1+\hat{r}_{2}-\hat{r}_{1}\cos
2\sigma \right) \hat{\Delta}_{12}\left( \hat{r}_{1},\hat{r}_{2}\right)
+3\left( \hat{r}_{2}-\hat{r}_{1}\cos 2\sigma \right) \left[ \left( \hat{r}%
_{2}-\hat{r}_{1}\right) +2\hat{r}_{1}\sin ^{2}\sigma \right] \right) d\sigma
}{\mathring{\Delta}_{12}^{5/2}\left( \hat{r}_{1},\hat{r}_{2}\right) },\;
\notag \\
\partial _{x_{1}}\Psi _{1}\left( q\right) & =\partial _{x_{2}}\Psi
_{1}\left( q\right) =0,  \notag \\
\partial _{s_{1}}\Psi _{2}\left( q\right) & =\frac{1}{\hat{r}_{1}}%
\int\nolimits_{0}^{\pi /2}\frac{\left( \left( 1+\hat{r}_{1}-\hat{r}_{2}\cos
2\sigma \right) \hat{\Delta}_{12}\left( \hat{r}_{1},\hat{r}_{2}\right)
+3\left( \hat{r}_{1}-\hat{r}_{2}\cos 2\sigma \right) \left[ \left( \hat{r}%
_{1}-\hat{r}_{2}\right) +2\hat{r}_{2}\sin ^{2}\sigma \right] \right) d\sigma
}{\mathring{\Delta}_{12}^{5/2}\left( \hat{r}_{1},\hat{r}_{2}\right) },\;
\notag \\
\partial _{x_{1}}\Psi _{2}\left( q\right) & =\partial _{x_{2}}\Psi
_{2}\left( q\right) =0,  \notag \\
\partial _{s_{2}}\Psi _{2}\left( q\right) & =-\frac{1}{2\hat{r}_{2}}\left\{
\frac{1}{2\hat{r}_{2}^{2}}\left[ \log \left( \chi \hat{r}_{2}\right) -\left(
1+\gamma \right) \right] -\right.  \notag \\
& \quad \left. 2\hat{r}_{1}\int\nolimits_{0}^{\pi /2}\frac{\left( \hat{\Delta%
}_{12}\left( \hat{r}_{1},\hat{r}_{2}\right) \cos 2\sigma +3\left( \hat{r}%
_{1}-\hat{r}_{2}\cos 2\sigma \right) \left[ \left( \hat{r}_{2}-\hat{r}%
_{1}\right) +2\hat{r}_{1}\sin ^{2}\sigma \right] \right) d\sigma }{\mathring{%
\Delta}_{12}^{5/2}\left( \hat{r}_{1},\hat{r}_{2}\right) }\right\} .
\label{e47}
\end{align}

To simplify our analysis of the spectrum of $\mathfrak{X}_{12}^{\prime
}\left( q\right) $ at each of the fixed points, we introduce the notation%
\begin{align*}
q_{\mathrm{I}}& =q_{\mathrm{I}}\left( \xi \right) :=\left( \hat{s}_{1(%
\mathrm{I})},\hat{s}_{2(\mathrm{I})},\xi ,\xi \right) ,\;q_{\mathrm{II}}=q_{%
\mathrm{II}}\left( \xi \right) :=\left( \hat{s}_{1(\mathrm{II})},\hat{s}_{2(%
\mathrm{II})},\xi ,\xi \right) , \\
q_{\mathrm{III}}& =q_{\mathrm{III}}\left( \xi \right) :=\left( \hat{s}_{1(%
\mathrm{III})},\hat{s}_{2(\mathrm{III})},\xi ,\xi \right) ,\;\mathrm{and}%
\;q_{\mathrm{IV}}=q_{\mathrm{IV}}\left( \xi \right) :=\left( \hat{s}_{1(%
\mathrm{IV})},\hat{s}_{2(\mathrm{IV})},\xi ,\xi \right) .
\end{align*}%
It is a routine - albeit tedious - matter to show that all four of the
derivatives $\mathfrak{X}_{12}^{\prime }\left( q_{I}\right) ,\ldots ,$ $%
\mathfrak{X}_{12}^{\prime }\left( q_{\mathrm{IV}}\right) $ have a pair of
zero eigenvalues, which is consistent with the fact that none of the
stationary points is isolated and all of the derivatives are symplectic
matrices, and a purely imaginary conjugate pair: $\lambda _{3(\mathrm{I}%
)}=i\nu _{\mathrm{I}}$, $\lambda _{4(\mathrm{I})}=\bar{\lambda}_{3(\mathrm{I}%
)}=-i\nu _{\mathrm{I}}$; $\lambda _{3(\mathrm{II})}=i\nu _{\mathrm{II}}$, $%
\lambda _{4(\mathrm{II})}=\bar{\lambda}_{3(\mathrm{II})}=-i\nu _{\mathrm{II}%
} $; $\lambda _{3(\mathrm{III})}=i\nu _{\mathrm{III}}$, $\lambda _{4(\mathrm{%
III})}=\bar{\lambda}_{3(\mathrm{III})}=-i\nu _{\mathrm{III}}$; and $\lambda
_{3(\mathrm{IV})}=i\nu _{\mathrm{IV}}$, $\lambda _{4(\mathrm{IV})}=\bar{%
\lambda}_{3(\mathrm{IV})}=-i\nu _{\mathrm{IV}}$, respectively, with $\nu _{%
\mathrm{I}},\nu _{\mathrm{II}},\nu _{\mathrm{III}},\nu _{\mathrm{IV}}>0$.

\section{Analysis of the Ring Kinematics}

With the ring dynamics disposed of in the preceding section, it is now quite
simple to resolve the ring kinematics. We begin our analysis from the
passive third ring approach as described in Subsection 3.1. Our calculations
are rendered simpler owing to the fact that the first four equations of (20)
do not depend on the last two equations, so the analysis of stationary
points for the ring dynamics can be used directly for the kinematics.

\subsection{Stationary points in the passive third ring approach}

First we rewrite (20) incorporating the assumptions of Section 4 and
introducing some additional notation:

\begin{align}
\dot{s}_{1}& =\Phi _{1}\left( r_{1},r_{2},x_{1},x_{2};\alpha ,\kappa \right)
:=4\kappa r_{1}r_{2}\left( x_{1}-x_{2}\right) \int\nolimits_{0}^{\pi /2}%
\frac{\cos 2\sigma d\sigma }{\Delta _{12}^{3/2}},  \notag \\
\dot{s}_{2}& =\Phi _{2}\left( r_{1},r_{2},x_{1},x_{2};\alpha ,\kappa \right)
:=4r_{1}r_{2}\left( x_{2}-x_{1}\right) \int\nolimits_{0}^{\pi /2}\frac{\cos
2\sigma d\sigma }{\Delta _{12}^{3/2}},  \notag \\
\dot{x}_{1}& =\Psi _{1}\left( r_{1},r_{2},x_{1},x_{2};\alpha ,\kappa \right)
:=-\alpha +\frac{1}{2r_{1}}\left[ \log \left( \chi r_{1}\right) -\gamma %
\right] +2\kappa r_{2}\int\nolimits_{0}^{\pi /2}\frac{\left( r_{2}-r_{1}\cos
2\sigma \right) d\sigma }{\Delta _{12}^{3/2}},  \notag \\
\dot{x}_{2}& =\Psi _{2}\left( r_{1},r_{2},x_{1},x_{2};\alpha ,\kappa \right)
:=-\alpha +\frac{\kappa }{2r_{2}}\left[ \log \left( \chi r_{2}\right)
-\gamma \right] +2r_{1}\int\nolimits_{0}^{\pi /2}\frac{\left(
r_{1}-r_{2}\cos 2\sigma \right) d\sigma }{\Delta _{12}^{3/2}},  \notag \\
\dot{s}& =\Phi \left( r_{1},r_{2},r,x_{1},x_{2},x;\alpha ,\kappa \right)
:=4r\sum\limits_{k=1}^{2}\kappa _{k}r_{k}\left( x-x_{k}\right)
\int\nolimits_{0}^{\pi /2}\frac{\cos 2\sigma d\sigma }{\Delta _{k}^{3/2}},
\notag \\
\dot{x}& =\Psi \left( r_{1},r_{2},r,x_{1},x_{2},x;\alpha ,\kappa \right)
:=-\alpha +2\sum\limits_{k=1}^{2}\kappa _{k}r_{k}\int\nolimits_{0}^{\pi /2}%
\frac{\left( r_{k}-r\cos 2\sigma \right) d\sigma }{\Delta _{k}^{3/2}},
\label{e48}
\end{align}%
where $\kappa _{1}=1<\kappa :=\kappa _{2}$. To find the set of stationary
points, we set all of these equations equal to zero. But we have already
taken care of the first four equations; namely, we found the line $\mathfrak{%
L}$ defined by (45) comprised of a 1-parameter infinity of (stationary)
points for which $\Phi _{1}=\Phi _{2}=\Psi _{1}=\Psi _{2}=0$. This leaves us
to solve
\begin{align}
\Phi & :=4r\left( x-\xi \right) \left[ \hat{r}_{1}\int\nolimits_{0}^{\pi /2}%
\frac{\cos 2\sigma d\sigma }{\hat{\Delta}_{1}^{3/2}}+\kappa \hat{r}%
_{2}\int\nolimits_{0}^{\pi /2}\frac{\cos 2\sigma d\sigma }{\hat{\Delta}%
_{2}^{3/2}}\right] =0,  \notag \\
\Psi & =-\alpha +2\left[ \hat{r}_{1}\int\nolimits_{0}^{\pi /2}\frac{\left(
\hat{r}_{1}-r\cos 2\sigma \right) d\sigma }{\hat{\Delta}_{1}^{3/2}}+\kappa
\hat{r}_{2}\int\nolimits_{0}^{\pi /2}\frac{\left( \hat{r}_{2}-r\cos 2\sigma
\right) d\sigma }{\hat{\Delta}_{2}^{3/2}}\right] =0,  \label{e49}
\end{align}%
where

\begin{equation}
\hat{\Delta}_{k}:=\left( r-\hat{r}_{k}\right) ^{2}+\left( x-\xi \right)
^{2}+4r\hat{r}_{k}\sin ^{2}\sigma ,\quad (1\leq k\leq 2)  \label{e50}
\end{equation}%
for all choices of pairs of $\hat{r}_{1}$ and $\hat{r}_{2}$ included in $%
\mathfrak{L}$; namely, $\hat{r}_{1(\mathrm{I})}$,\ldots , $\hat{r}_{1(%
\mathrm{IV})}$ and $\hat{r}_{2(\mathrm{I})}$,\ldots , $\hat{r}_{2(\mathrm{IV}%
)}$, respectively. It is therefore helpful to identify the following cases
for the coaxial ring fixed point types:

\begin{center}
\begin{tabular}{ll}
Type I. & $r_{1}=\hat{r}_{1}=\hat{r}_{1(\mathrm{I})}$, $r_{2}=\hat{r}_{2}=%
\hat{r}_{2(\mathrm{I})}$ \\
Type II. & $r_{1}=\hat{r}_{1}=\hat{r}_{1(\mathrm{II})}$, $r_{2}=\hat{r}_{2}=%
\hat{r}_{2(\mathrm{II})}$ \\
Type III. & $r_{1}=\hat{r}_{1}=\hat{r}_{1(\mathrm{III})}$, $r_{2}=\hat{r}%
_{2}=\hat{r}_{2(\mathrm{III})}$ \\
Type IV. & $r_{1}=\hat{r}_{1}=\hat{r}_{1(\mathrm{IV})}$, $r_{2}=\hat{r}_{2}=%
\hat{r}_{2(\mathrm{IV})}$%
\end{tabular}
\end{center}

As the bracketed term in the first equation of (49) is positive, the only
choices for zeros are $r=0$ and $x=\xi $. If $r=0$, and we consider Type I,
the second equation of (49) reduces to the form%
\begin{equation}
\Psi =-\alpha +\pi \left\{ \frac{\hat{r}_{1(\mathrm{I})}^{2}}{\left[ \hat{r}%
_{1(\mathrm{I})}^{2}+\left( x-\xi \right) ^{2}\right] ^{3/2}}+\frac{\kappa
\hat{r}_{2(\mathrm{I})}^{2}}{\left[ \hat{r}_{2(\mathrm{I})}^{2}+\left( x-\xi
\right) ^{2}\right] ^{3/2}}\right\} =0,  \label{e51}
\end{equation}%
which immediately yields a pair of solutions that are symmetric with respect
to any chosen value of $x=\xi $; namely,%
\begin{equation}
x_{\mathrm{I}}^{(\pm )}:=\xi \pm \eta _{\mathrm{I}},  \label{e52}
\end{equation}%
where $\eta _{\mathrm{I}}>0$. We note that $x_{\mathrm{I}}^{(+)}$ and $x_{%
\mathrm{I}}^{(-)}$ correspond for Type I, respectively, to the stagnation
points on the leading and trailing edges of the vortex breakdown bubble (on
the axis of symmetry). Analogously, we find points $x_{\mathrm{II}}^{(\pm
)}:=\xi \pm \eta _{\mathrm{II}}$, $x_{\mathrm{III}}^{(\pm )}:=\xi \pm \eta _{%
\mathrm{III}}$ and $x_{\mathrm{IV}}^{(\pm )}:=\xi \pm \eta _{\mathrm{IV}}$
for Types II, III and IV, respectively.

If we choose $x=\xi $ for Type I (for $\alpha >1$) , the second equation of
(49) becomes%
\begin{equation}
\Psi =-\alpha +2\Xi \left( r;\hat{r}_{1(\mathrm{I})},\hat{r}_{2(\mathrm{I}%
)},\kappa \right) =0,  \label{e53}
\end{equation}%
where%
\begin{equation}
\Xi \left( r;\hat{r}_{1(\mathrm{I})},\hat{r}_{2(\mathrm{I})},\kappa \right)
:=\int\nolimits_{0}^{\pi /2}\left\{ \frac{\hat{r}_{1(\mathrm{I})}\left( \hat{%
r}_{1(\mathrm{I})}-r\cos 2\sigma \right) }{\left[ \left( r-\hat{r}_{1(%
\mathrm{I})}\right) ^{2}+4r\hat{r}_{1(\mathrm{I})}\sin ^{2}\sigma \right]
^{3/2}}+\frac{\kappa \hat{r}_{2(\mathrm{I})}\left( \hat{r}_{2(\mathrm{I}%
)}-r\cos 2\sigma \right) }{\left[ \left( r-\hat{r}_{2(\mathrm{I})}\right)
^{2}+4r\hat{r}_{2(\mathrm{I})}\sin ^{2}\sigma \right] ^{3/2}}\right\}
d\sigma .  \label{e54}
\end{equation}%
Then a routine - but rather laborious - analysis of the properties of $\Xi $
shows that it is a smooth, nondecreasing function of $r$ for $\hat{r}_{1(%
\mathrm{I})}<r<\hat{r}_{2(\mathrm{I})}$ such that $\lim_{r\downarrow \hat{r}%
_{1(\mathrm{I})}}\Xi =-\infty $, $\lim_{r\uparrow \hat{r}_{2(\mathrm{I}%
)}}\Xi =\infty $, and the derivative with respect to $r$ increases rapidly
from nearly zero to $+\infty $ very close to $\hat{r}_{2(\mathrm{I})}$. It
also is easy to verify that $\Xi $ is nonvanishing when $r\notin \lbrack
\hat{r}_{1(\mathrm{I})},\hat{r}_{2(\mathrm{I})}]$. Therefore, in light of
the definitions of $\hat{r}_{1(\mathrm{I})}$ and $\hat{r}_{2(\mathrm{I})}$
as a solution pair of (39) with $\hat{r}_{1(\mathrm{I})}<\hat{r}_{2(\mathrm{I%
})}\simeq 4\hat{r}_{1(\mathrm{I})}$, we infer that (54) has a unique
solution $\hat{r}_{\mathrm{I}}$, which is increasingly near to $\hat{r}_{1(%
\mathrm{I})}$ as $\alpha $ increases; in particular,
\begin{equation*}
0<\hat{r}_{1(\mathrm{I})}<\hat{r}_{\mathrm{I}}<\hat{r}_{2(\mathrm{I})}\;%
\mathrm{and}\;\hat{r}_{\mathrm{I}}\downarrow \hat{r}_{1(\mathrm{I})}\;%
\mathrm{as}\;\alpha \uparrow \infty .
\end{equation*}%
Observe that one should also include the singular stationary points $\left(
s,x\right) =\left( \hat{s}_{1(\mathrm{I})},\xi \right) $ and $\left(
s,x\right) =\left( \hat{s}_{2(\mathrm{I})},\xi \right) $, which are not
actually included in the phase space - but are in fact stationary with
respect to the kinematics - and can readily be shown to behave like
(singular) centers surrounded locally by periodic orbits. It is
straightforward to verify the analogous behavior for solutions of (54) for
Types II, III, and IV; in particular, one finds the following: For Type II
there is a unique solution $\hat{r}_{\mathrm{II}}$ such that $\hat{r}_{2(%
\mathrm{II})}<\hat{r}_{\mathrm{II}}<\hat{r}_{1(\mathrm{II})}$; there is a
unique solution $\hat{r}_{\mathrm{III}}$ with $\hat{r}_{1(\mathrm{III})}<%
\hat{r}_{\mathrm{III}}<\hat{r}_{2(\mathrm{III})}$ for Type III; and for Type
IV there is a unique solution $\hat{r}_{\mathrm{IV}}$ such that $\hat{r}_{1(%
\mathrm{IV})}<\hat{r}_{\mathrm{IV}}<\hat{r}_{2(\mathrm{IV})}$.

To summarize our findings here, we have found that the complete set $%
\mathcal{F}$ of stationary points of (20), including the singular stationary
points, in the 6-dimensional, $s_{1},s_{2},s,x_{1},x_{2},x$-(phase) space $%
\boldsymbol{X}_{DK}$, is the following 1-parameter family:%
\begin{align}
\mathcal{F}& =\mathcal{F}\left( \xi \right) :=\left\{ \left( \hat{s}_{1(%
\mathrm{I})}=\hat{r}_{1(\mathrm{I})}^{2},\hat{s}_{2(\mathrm{I})}=\hat{r}_{2(%
\mathrm{I})}^{2},0,\xi ,\xi ,\xi \pm \eta _{\mathrm{I}}\right) ,\left( \hat{s%
}_{1(\mathrm{II})}=\hat{r}_{1(\mathrm{II})}^{2},\hat{s}_{2(\mathrm{II})}=%
\hat{r}_{2(\mathrm{II})}^{2},0,\xi ,\xi ,\xi \pm \eta _{\mathrm{II}}\right)
,\right.  \notag \\
& \qquad \qquad \qquad \left. \left( \hat{s}_{1(\mathrm{III})}=\hat{r}_{1(%
\mathrm{III})}^{2},\hat{s}_{2(\mathrm{III})}=\hat{r}_{2(\mathrm{III}%
)}^{2},0,\xi ,\xi ,\xi \pm \eta _{\mathrm{III}}\right) ,\left( \hat{s}_{1(%
\mathrm{IV})}=\hat{r}_{1(\mathrm{IV})}^{2},\hat{s}_{2(\mathrm{IV})}=\hat{r}%
_{2(\mathrm{IV})}^{2},0,\xi ,\xi ,\xi \pm \eta _{\mathrm{IV}}\right) ,\right.
\notag \\
& \qquad \qquad \qquad \quad \left. \left( \hat{s}_{1(\mathrm{I})},\hat{s}%
_{2(\mathrm{I})},\hat{s}_{\mathrm{I}}=\hat{r}_{\mathrm{I}}^{2},\xi ,\xi ,\xi
\right) ,\left( \hat{s}_{1(\mathrm{II})},\hat{s}_{2(\mathrm{II})},\hat{s}_{%
\mathrm{II}}=\hat{r}_{\mathrm{II}}^{2},\xi ,\xi ,\xi \right) ,\right.  \notag
\\
& \qquad \qquad \qquad \quad \quad \left. \left( \hat{s}_{1(\mathrm{III})},%
\hat{s}_{2(\mathrm{III})},\hat{s}_{\mathrm{III}}=\hat{r}_{\mathrm{III}%
}^{2},\xi ,\xi ,\xi \right) ,\left( \hat{s}_{1(\mathrm{IV})},\hat{s}_{2(%
\mathrm{IV})},\hat{s}_{\mathrm{IV}}=\hat{r}_{\mathrm{IV}}^{2},\xi ,\xi ,\xi
\right) \right.  \notag \\
& \qquad \qquad \qquad \quad \quad \quad \left. \left( \hat{s}_{1(\mathrm{I}%
)},\hat{s}_{2(\mathrm{I})},\hat{s}_{1(\mathrm{I})},\xi ,\xi ,\xi \right)
,\left( \hat{s}_{1(\mathrm{I})},\hat{s}_{2(\mathrm{I})},\hat{s}_{2(\mathrm{I}%
)},\xi ,\xi ,\xi \right) ,\left( \hat{s}_{1(\mathrm{II})},\hat{s}_{2(\mathrm{%
II})},\hat{s}_{1(\mathrm{II})},\xi ,\xi ,\xi \right) ,\right.  \notag \\
& \qquad \qquad \qquad \quad \quad \qquad \left. \left( \hat{s}_{1(\mathrm{II%
})},\hat{s}_{2(\mathrm{II})},\hat{s}_{2(\mathrm{II})},\xi ,\xi ,\xi \right)
,\left( \hat{s}_{1(\mathrm{III})},\hat{s}_{2(\mathrm{III})},\hat{s}_{1(%
\mathrm{III})},\xi ,\xi ,\xi \right) ,\left( \hat{s}_{1(\mathrm{III})},\hat{s%
}_{2(\mathrm{III})},\hat{s}_{2(\mathrm{III})},\xi ,\xi ,\xi \right) ,\right.
\notag \\
& \qquad \qquad \qquad \quad \quad \qquad \qquad \qquad \qquad \left. \left(
\hat{s}_{1(\mathrm{IV})},\hat{s}_{2(\mathrm{IV})},\hat{s}_{1(\mathrm{IV}%
)},\xi ,\xi ,\xi \right) ,\left( \hat{s}_{1(\mathrm{IV})},\hat{s}_{2(\mathrm{%
IV})},\hat{s}_{2(\mathrm{IV})},\xi ,\xi ,\xi \right) :\xi \in \mathbb{R}%
\right\}  \label{e55}
\end{align}

\subsection{Characterization of stationary points}

Owing to our analysis in Subsection 4.2, and the nature of the system of
differential equations under investigation, it suffices to determine the
type of the stationary points with respect to one of the invariant
half-planes of the following forms, depending upon which stationary points
are chosen for the pair of coaxial rings as specified by the four cases
above:%
\begin{align}
\mathcal{S}_{\mathrm{I}}\left( \xi \right) & :=\left\{ P=\left( \hat{s}_{1(%
\mathrm{I})},\hat{s}_{2(\mathrm{I})},s,\xi ,\xi ,x\right) :s\geq 0,\;x\in
\mathbb{R}\right\} ,  \notag \\
\mathcal{S}_{\mathrm{II}}\left( \xi \right) & :=\left\{ P=\left( \hat{s}_{1(%
\mathrm{II})},\hat{s}_{2(\mathrm{II})},s,\xi ,\xi ,x\right) :s\geq 0,\;x\in
\mathbb{R}\right\} ,  \notag \\
\mathcal{S}_{\mathrm{III}}\left( \xi \right) & :=\left\{ P=\left( \hat{s}_{1(%
\mathrm{III})},\hat{s}_{2(\mathrm{III})},s,\xi ,\xi ,x\right) :s\geq
0,\;x\in \mathbb{R}\right\} ,  \notag \\
\mathcal{S}_{\mathrm{IV}}\left( \xi \right) & :=\left\{ P=\left( \hat{s}_{1(%
\mathrm{IV})},\hat{s}_{2(\mathrm{IV})},s,\xi ,\xi ,x\right) :s\geq 0,\;x\in
\mathbb{R}\right\}  \label{e56}
\end{align}%
To do the linear analysis, we have to compute the derivatives of the vector
field in these half-planes at the stationary points $P\in \mathcal{S}\left(
\xi \right) $; namely,%
\begin{equation}
\mathcal{X}^{\prime }(P):=\left(
\begin{array}{cc}
\partial _{s}\Phi \left( P\right) & \partial _{x}\Phi \left( P\right) \\
\partial _{s}\Psi \left( P\right) & \partial _{x}\Psi \left( P\right)%
\end{array}%
\right) .  \label{e57}
\end{equation}%
We shall do this in rather complete detail for Type I (with $\alpha >1$),
then in somewhat less detail for Type II, and then just briefly summarize
the analysis for Types III and IV. Toward this end, we readily obtain the
following formulas:%
\begin{align}
\partial _{s}\Phi \left( P\right) & =\frac{2\left( x-\xi \right) }{r}%
\sum\limits_{k=1}^{2}\int\nolimits_{0}^{\pi /2}\frac{\kappa _{k}\hat{r}_{k}%
\left[ \left( x-\xi \right) ^{2}-\left( r-\hat{r}_{k}\right) \left( 2r+\hat{r%
}_{k}\right) -2r\hat{r}_{k}\sin ^{2}\sigma \right] \cos 2\sigma d\sigma }{%
\left[ \left( r-\hat{r}_{k}\right) ^{2}+\left( x-\xi \right) ^{2}+4r\hat{r}%
_{k}\sin ^{2}\sigma \right] ^{5/2}},  \notag \\
\partial _{x}\Phi \left( P\right) &
=4r\sum\limits_{k=1}^{2}\int\nolimits_{0}^{\pi /2}\frac{\kappa _{k}\hat{r}%
_{k}\left[ \left( r-\hat{r}_{k}\right) ^{2}+4r\hat{r}_{k}\sin ^{2}\sigma
-2\left( x-\xi \right) ^{2}\right] \cos 2\sigma d\sigma }{\left[ \left( r-%
\hat{r}_{k}\right) ^{2}+\left( x-\xi \right) ^{2}+4r\hat{r}_{k}\sin
^{2}\sigma \right] ^{5/2}},  \notag \\
\partial _{s}\Psi \left( P\right) & =\frac{1}{r}\sum\limits_{k=1}^{2}\int%
\nolimits_{0}^{\pi /2}\frac{\kappa _{k}\hat{r}_{k}^{2}\left\{ \left[ \left(
2r^{2}-r\hat{r}_{k}-\hat{r}_{k}^{2}\right) -\left( x-\xi \right) ^{2}+2r\hat{%
r}_{k}\sin ^{2}\sigma \right] \cos 2\sigma -6r\hat{r}_{k}\sin ^{2}\sigma
\right\} d\sigma }{\left[ \left( r-\hat{r}_{k}\right) ^{2}+\left( x-\xi
\right) ^{2}+4r\hat{r}_{k}\sin ^{2}\sigma \right] ^{5/2}},  \notag \\
\partial _{x}\Psi \left( P\right) & =-6\left( x-\xi \right)
\sum\limits_{k=1}^{2}\int\nolimits_{0}^{\pi /2}\frac{\kappa _{k}\hat{r}%
_{k}\left( \hat{r}_{k}-r\cos 2\sigma \right) d\sigma }{\left[ \left( r-\hat{r%
}_{k}\right) ^{2}+\left( x-\xi \right) ^{2}+4r\hat{r}_{k}\sin ^{2}\sigma %
\right] ^{5/2}}.  \label{e58}
\end{align}

Now we fix a particular value of $\xi $, which we denote as $\hat{\xi}$, and
proceed to our line of analysis of the stationary points in $\mathcal{S}_{%
\mathrm{I}}\left( \hat{\xi}\right) $. First, we investigate the derivative
at the fixed point pair $P=p_{\mathrm{I}_{\pm }}:=\left( \hat{s}_{1(\mathrm{I%
})}=\hat{r}_{1(\mathrm{I})}^{2},\hat{s}_{2(\mathrm{I})}=\hat{r}_{2(\mathrm{I}%
)}^{2},0,\hat{\xi},\hat{\xi},\hat{\xi}\pm \eta _{\mathrm{I}}\right) $. We
compute from (58) that%
\begin{equation}
\mathcal{X}^{\prime }(p_{\pm })=\left(
\begin{array}{cc}
\partial _{s}\Phi \left( p_{\mathrm{I}_{\pm }}\right) & 0 \\
\partial _{s}\Psi \left( p_{\mathrm{I}_{\pm }}\right) & \partial _{x}\Psi
\left( p_{\mathrm{I}_{\pm }}\right)%
\end{array}%
\right) ,  \label{e59}
\end{equation}%
where%
\begin{align}
\partial _{s}\Phi \left( p_{\mathrm{I}_{\pm }}\right) & =\pm \frac{\pi \eta
_{\mathrm{I}}}{2}\left[ \frac{\hat{r}_{1(\mathrm{I})}^{2}}{\left( \hat{r}_{1(%
\mathrm{I})}^{2}+\eta _{\mathrm{I}}^{2}\right) ^{5/2}}+\frac{\kappa \hat{r}%
_{2(\mathrm{I})}^{2}}{\left( \hat{r}_{2(\mathrm{I})}^{2}+\eta _{\mathrm{I}%
}^{2}\right) ^{5/2}}\right] ,  \notag \\
\partial _{s}\Psi \left( p_{\mathrm{I}_{\pm }}\right) & =-\pi \left[ \frac{%
\hat{r}_{1(\mathrm{I})}^{3}}{\left( \hat{r}_{1(\mathrm{I})}^{2}+\eta _{%
\mathrm{I}}^{2}\right) ^{5/2}}+\frac{\kappa \hat{r}_{2(\mathrm{I})}^{3}}{%
\left( \hat{r}_{2(\mathrm{I})}^{2}+\eta _{\mathrm{I}}^{2}\right) ^{5/2}}%
\right] ,  \label{e60} \\
\partial _{x}\Psi \left( p_{\mathrm{I}_{\pm }}\right) & =\mp 3\pi \eta _{%
\mathrm{I}}\left[ \frac{\hat{r}_{1(\mathrm{I})}^{2}}{\left( \hat{r}_{1(%
\mathrm{I})}^{2}+\eta _{\mathrm{I}}^{2}\right) ^{5/2}}+\frac{\kappa \hat{r}%
_{2(\mathrm{I})}^{2}}{\left( \hat{r}_{2(\mathrm{I})}^{2}+\eta _{\mathrm{I}%
}^{2}\right) ^{5/2}}\right] ,  \notag
\end{align}%
where it should be noted that the partial derivatives with respect to $s$
need to be calculated in a limiting sense as $r\downarrow 0$. Whence, we
immediately conclude that both $p_{\mathrm{I}_{+}}$ and $p_{\mathrm{I}_{-}}$
are saddle points, with stable and unstable manifolds $W^{s}$ and $W^{u}$ in
$\mathcal{S}_{\mathrm{I}}\left( \xi \right) $, respectively, satisfying the
following properties :%
\begin{align}
W^{s}\left( p_{\mathrm{I}_{+}}\right) & =\left( x_{\mathrm{I}}^{(-)},\infty
\right) ,  \notag \\
W^{u}\left( p_{_{\mathrm{I}_{-}}}\right) & =\left( -\infty ,x_{\mathrm{I}%
}^{(+)}\right) ,  \label{e61}
\end{align}%
and $W^{u}\left( p_{\mathrm{I}_{+}}\right) $ at $p_{\mathrm{I}_{+}}$ is
tangent to the line with nonpositive slope given as%
\begin{equation}
L_{\mathrm{I}_{+}}:\partial _{s}\Psi \left( p_{\mathrm{I}_{+}}\right) s+%
\left[ \partial _{x}\Psi \left( p_{\mathrm{I}_{+}}\right) -\partial _{s}\Phi
\left( p_{\mathrm{I}_{+}}\right) \right] x=0,  \label{e62}
\end{equation}%
while $W^{s}\left( p_{\mathrm{I}_{-}}\right) $ at $p_{\mathrm{I}_{-}}$ is
tangent to the line of nonnegative slope of the form%
\begin{equation}
L_{\mathrm{I}_{-}}:\partial _{s}\Psi \left( p_{\mathrm{I}_{-}}\right) s+%
\left[ \partial _{x}\Psi \left( p_{\mathrm{I}_{-}}\right) -\partial _{s}\Phi
\left( p_{\mathrm{I}_{-}}\right) \right] x=0.  \label{e63}
\end{equation}%
Moreover, it is easy to see from the form of the dynamical system (48) that $%
W^{u}\left( p_{\mathrm{I}_{+}}\right) \smallsetminus \{p_{\mathrm{I}%
_{+}}\}=W^{s}\left( p_{\mathrm{I}_{-}}\right) \smallsetminus \{p_{\mathrm{I}%
_{-}}\}:=\mathcal{C}_{\mathrm{I}}:=\mathcal{C}_{\mathrm{I}}\left( p_{\mathrm{%
I}_{+}},p_{\mathrm{I}_{-}}\right) $ is a smooth, convex curve joining $p_{%
\mathrm{I}_{-}}$ to $p_{\mathrm{I}_{+}}$, is symmetric with respect to the
line $x=\hat{\xi}$, and is (asymptotically) tangent to $L_{\mathrm{I}_{+}}$
and $L_{\mathrm{I}_{-}}$ at $p_{\mathrm{I}_{+}}$ and $p_{\mathrm{I}_{-}}$,
respectively (see Fig. 1). Thus, $\mathcal{C}_{\mathrm{I}}$ together with
the line segment $\{\left( 0,x\right) :x_{\mathrm{I}}^{(-)}\leq x\leq x_{%
\mathrm{I}}^{(+)}\}$ comprises a heteroclinic cycle that we denote as $%
\mathcal{Z}_{\mathrm{I}}$, which is susceptible to perturbations that can
generate chaotic dynamical regimes, such as described in \cite{bb}, \cite%
{dbiutam}, \cite{BB}, \cite{guho}, \cite{kathas}, \cite{malwig}, \cite{MSWI}%
, \cite{sotir2}, and \cite{wigbook}.

Next we consider the stationary points in $\mathcal{F}$ of the form $%
Q=\left(  \hat{s}_{1},\hat{s}_{2},s,\hat{\xi},\hat{\xi},\hat{\xi}\right) $,
where we have chosen a specific value $\hat{\xi}$ of $\xi $ for our
computations. As we are considering Case I, $\hat{s}_{1}=\hat{s}_{1(\mathrm{%
I })}$ and $\hat{s}_{2}=\hat{s}_{2(\mathrm{I})}$. Again we use formulas (57)
and (58) to calculate the (linear) type of the stationary point via $%
\mathcal{X}^{\prime }(Q)$. The following results may readily be verified:
\begin{align}
\partial _{s}\Phi \left( Q\right) & =\partial _{x}\Psi \left( Q\right) =0,
\notag \\
\partial _{x}\Phi \left( Q\right) &
=4r\sum\limits_{k=1}^{2}\int\nolimits_{0}^{\pi /2}\frac{\kappa _{k}\hat{r}%
_{k(\mathrm{I})}\cos 2\sigma d\sigma }{\left[ \left( r-\hat{r}_{k(\mathrm{I}%
)}\right) ^{2}+4r\hat{r}_{k(\mathrm{I})}\sin ^{2}\sigma \right] ^{5/2}},
\label{e64} \\
\partial _{s}\Psi \left( Q\right) & =\frac{1}{r}\sum\limits_{k=1}^{2}\int%
\nolimits_{0}^{\pi /2}\frac{\kappa _{k}\hat{r}_{k(\mathrm{I})}^{2}\left[
\left( 2r^{2}-r\hat{r}_{k(\mathrm{I})}-\hat{r}_{k(\mathrm{I})}^{2}+2r\hat{r}%
_{k(\mathrm{I})}\sin ^{2}\sigma \right) \cos 2\sigma -6r\hat{r}_{k(\mathrm{I}%
)}\sin ^{2}\sigma \right] d\sigma }{\left[ \left( r-\hat{r}_{k(\mathrm{I}%
)}\right) ^{2}+4r\hat{r}_{k(\mathrm{I})}\sin ^{2}\sigma \right] ^{5/2}}.
\notag
\end{align}%
Observe that, as was to be expected, these formulas are independent of $\xi $%
.

It is easy to see that $\partial _{x}\Phi \left( Q\right) $ is positive for
all $r>0$ (and also all $\xi )$, and it can be shown by a routine but rather
lengthy analysis that $\partial _{s}\Psi \left( Q\right) $ is positive at $%
Q=q_{\mathrm{I}}:=\left( \hat{s}_{1(\mathrm{I})},\hat{s}_{2(\mathrm{I})},%
\hat{s}_{\mathrm{I}},\hat{\xi},\hat{\xi},\hat{\xi}\right) $. Accordingly the
eigenvalues of $\mathcal{X}^{\prime }(q_{\mathrm{I}})$ are $\pm \sqrt{%
\partial _{x}\Phi \left( q_{\mathrm{I}}\right) \partial _{s}\Psi \left( q_{%
\mathrm{I}}\right) }$, so $q_{\mathrm{I}}$ is a saddle point. As for the
singular stationary points $Q=q_{\mathrm{I}_{-}}:=\left( \hat{s}_{1(\mathrm{I%
})},\hat{s}_{2(\mathrm{I})},\hat{s}_{1(\mathrm{I})},\hat{\xi},\hat{\xi},\hat{%
\xi}\right) $ and $Q=q_{\mathrm{I}_{+}}:=\left( \hat{s}_{1(\mathrm{I})},\hat{%
s}_{2(\mathrm{I})},\hat{s}_{2(\mathrm{I})},\hat{\xi},\hat{\xi},\hat{\xi}%
\right) $, it is intuitively clear and follows from a straightforward
analysis of the (48), in which we employ the results obtained in Subsection
4.2, that these points behave, except right at the points, just like centers
of the Hamiltonian dynamical system.

We may now, after some additional routine analysis of the system, piece
together the results obtained on the stationary points and invariant curves
in order to produce a complete qualitative picture of the nature of the
phase portrait of (48) in any of the half-planes $\mathcal{S}_{\mathrm{I}%
}\left( \xi \right) $, which is illustrated in Fig. 1. The main details of
the phase planes includes the heteroclinic cycle $\mathcal{Z}_{\mathrm{I}}$
connecting the saddle points $p_{\mathrm{I}_{\pm }}$ on the $x$-axis, which
encloses the remaining three stationary points $q_{\mathrm{I}_{+}}$, $q_{%
\mathrm{I}_{-}}$, and $q_{\mathrm{I}}$ and all of their associated stable
and unstable manifolds. In particular, we see that
\begin{equation*}
\lbrack \{q_{\mathrm{I}}\}\cup W^{s}(q_{\mathrm{I}})\cup W^{u}(q_{\mathrm{I}%
})]\cap \{r\leq \hat{r}_{\mathrm{I}}\}
\end{equation*}%
forms a small homoclinic loop $\ell _{\mathrm{I}_{-}}$ enclosing $q_{\mathrm{%
I}_{-}}$, and
\begin{equation*}
\lbrack \{q_{\mathrm{I}}\}\cup W^{s}(q_{\mathrm{I}})\cup W^{u}(q_{\mathrm{I}%
})]\cap \{r\geq \hat{r}_{\mathrm{I}}\}
\end{equation*}%
forms a small homoclinic loop $\ell _{\mathrm{I}_{+}}$ enclosing $q_{\mathrm{%
I}_{+}}$. Moreover, it is not difficult to verify that
\begin{equation*}
\lbrack \{q_{\mathrm{I}}\}\cup W^{s}(q_{\mathrm{I}})\cup W^{u}(q_{\mathrm{I}%
})]\cap \{r>0\}
\end{equation*}%
forms a figure eight curve, which we denote as $\zeta _{\mathrm{I}}$,
enclosing the singular centers $q_{\mathrm{I}_{-}}$ and $q_{\mathrm{I}_{+}}$%
. Owing to the symmetry of the system with respect to the line $x=\hat{\xi}$%
, the homoclinic loops $\ell _{\mathrm{I}_{-}}$ and $\ell _{\mathrm{I}_{+}}$
and the curve $\zeta _{\mathrm{I}}$ are also symmetric under reflection in
this line. In addition, the interiors of $\ell _{\mathrm{I}_{-}}$ and $\ell
_{\mathrm{I}_{+}}$ are comprised of periodic orbits circling $q_{\mathrm{I}%
_{-}}$ and $q_{\mathrm{I}_{+}}$, respectively, while points close to $\zeta
_{\mathrm{I}}$ and in its exterior belong to periodic trajectories that
circle around it. Observe that this representation of the dynamics of (48)
is quite consistent with what one might expect for a stationary,
axisymmetric vortex breakdown flow of B-type. Thus Type I (when $\alpha >1$)
we have, with our two coaxial vortex ring model, produced a two-parameter ($%
\alpha $ and $\kappa $) family of axisymmetric B-type vortex breakdown
flows. These features are illustrated in Fig. 1.

For Type II, we find many similarities, but also some differences in the
phase (half-) plane $\mathcal{S}_{\mathrm{II}}\left( \hat{\xi}\right) $ in
comparison to the phase plane $\mathcal{S}_{\mathrm{I}}\left( \xi \right) $
for Type I. Once again we have saddle points, $p_{\mathrm{II}_{\pm }}$, on
the $x$-axis determined from (54) with the obvious modifications, where
their linear stability is determined via (60) \emph{mutatis mutandis}. We
also have the singular points, this time $q_{\mathrm{II}_{+}}:=\left( \hat{s}%
_{1(\mathrm{II})},\hat{s}_{2(\mathrm{II})},\hat{s}_{1(\mathrm{II})},\hat{\xi}%
,\hat{\xi},\hat{\xi}\right) $ and $q_{\mathrm{II}_{-}}:=\left( \hat{s}_{1(%
\mathrm{II})},\hat{s}_{2(\mathrm{II})},\hat{s}_{2(\mathrm{II})},\hat{\xi},%
\hat{\xi},\hat{\xi}\right) $, behaving like nonlinear centers. In addition,
there is another stationary point $q_{\mathrm{II}}:=\left( \hat{s}_{1(%
\mathrm{II})},\hat{s}_{2(\mathrm{II})},\hat{s}_{\mathrm{II}},\hat{\xi},\hat{%
\xi},\hat{\xi}\right) $, which lies between $q_{\mathrm{II}_{-}}$ and $q_{%
\mathrm{II}_{+}}$, and is a saddle point in virtue of (64) with the obvious
modifications. There is also a heteroclinic cycle $\mathcal{Z}_{\mathrm{II}}$
(describing the bubble)connecting $p_{\mathrm{II}_{-}}$ and $p_{\mathrm{II}%
_{+}}$, which is symmetric with respect to $x=\hat{\xi}$, and encloses the
stationary points $q_{\mathrm{II}_{-}}$ and $q_{\mathrm{II}_{+}}$. The
bubble shape produced for this type of fixed ring configuration tends to
look less like the one expected, as compared to Type I for $\alpha >1$
sufficiently large and Type IV as can be seen from Figs. 1, 2 and 4.

Type III is almost completely analogous to Type II, where the only changes
involve substituting $\hat{s}_{1(\mathrm{III})}$ for $\hat{s}_{1(\mathrm{II}%
)}$, $\hat{s}_{2(\mathrm{III})}$ for $\hat{s}_{2(\mathrm{II})}$, and $\hat{s}%
_{\mathrm{III}}$ for $\hat{s}_{\mathrm{II}}$. In particular, the phase
portrait in $\mathcal{S}_{\mathrm{III}}\left( \hat{\xi}\right) $ is
qualitatively the same as in Type II, consistent with the changes in the
stationary points that we just indicated. However, we see from Fig. 3 that
the bubble shape produced for this type of fixed coaxial vortex ring
configuration tends to look less like the physical B-type flow than any of
the other types, as can be confirmed by noting that the shape in Fig. 3
resembles a marble balanced on the axis of symmetry. Actually the shape for
Type IV, as depicted in Fig. 4, seems to be most like what one expects of a
physical vortex breakdown streamline pattern of B-type.

We note that in Figs. 2, 3 and 4, the small rectangles in the upper
right-hand corners are microscopic depictions of the much smaller scale
streamline topologies near the fixed coaxial vortex rings that are closest
to the axis of symmetry.



\begin{figure}[ht!]
\centering
\includegraphics[width=4in]{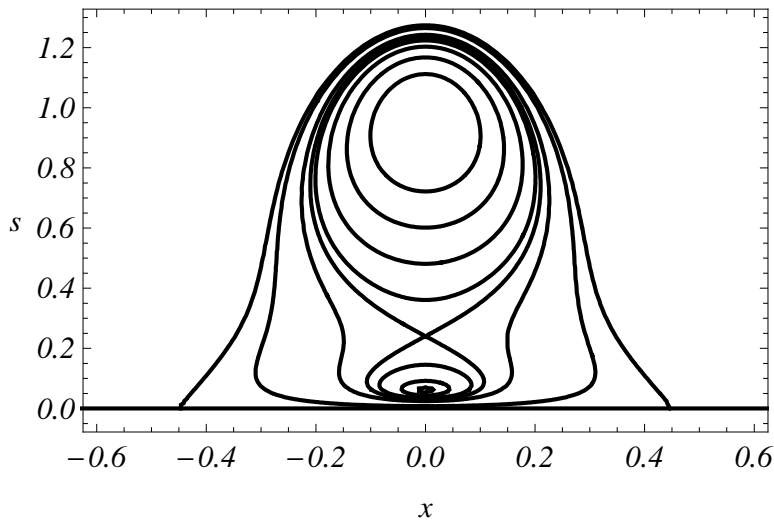}
\caption{Streamlines for Type I, $\protect\alpha=5,\; \protect\kappa=1.5$,
and $\hat{s}_1\simeq 0.06,\,\hat{s}_2\simeq 0.94,\, \hat{s} \simeq 0.24, \,
x_{\pm} \simeq\pm 0.45$.}
\label{streamlinecase1}
\end{figure}

\begin{figure}[ht!]
\centering
\includegraphics[width=4in]{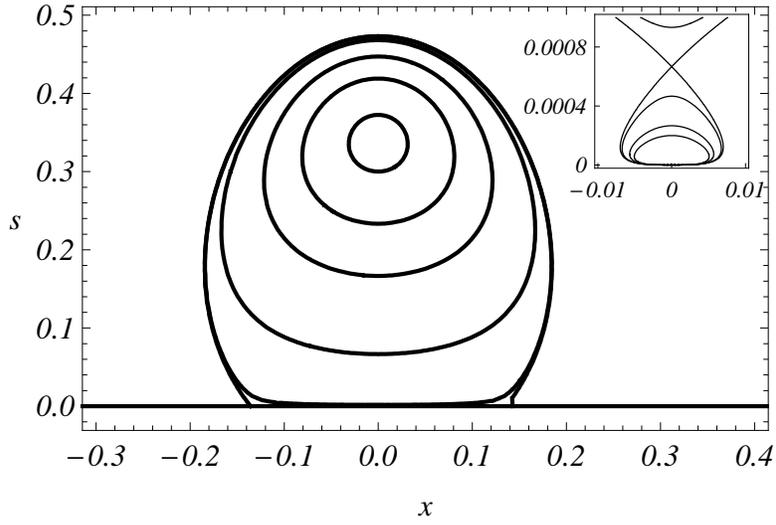}
\caption{Streamlines for Type II, $\protect\alpha=5,\; \protect\kappa=1.5$,
and $\hat{s}_1\simeq 0.33,\,\hat{s}_2\simeq 3\times 10^{-6},\, \hat{s}
\simeq 0.0006, \, x_{\pm} \simeq \pm 0.13$.}
\label{streamlinecase2}
\end{figure}

\begin{figure}[ht!]
\centering
\includegraphics[width=4in]{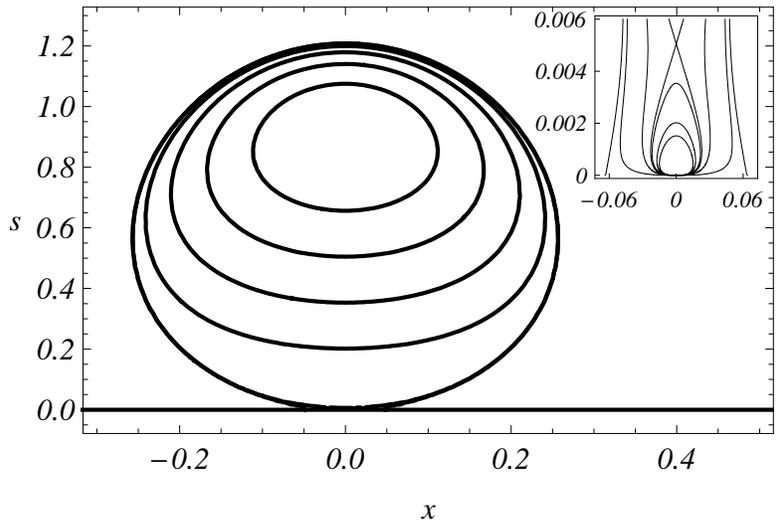}
\caption{Streamlines for Type III, $\protect\alpha=5,\; \protect\kappa=1.5$,
and $\hat{s}_1\simeq 3.2\times 10^{-6},\,\hat{s}_2\simeq 0.89,\, \hat{s}
\simeq 0.005 ,\, x_{\pm} \simeq \pm 0.063$.}
\label{streamlinecase3}
\end{figure}

\begin{figure}[ht!]
\centering
\includegraphics[width=4in]{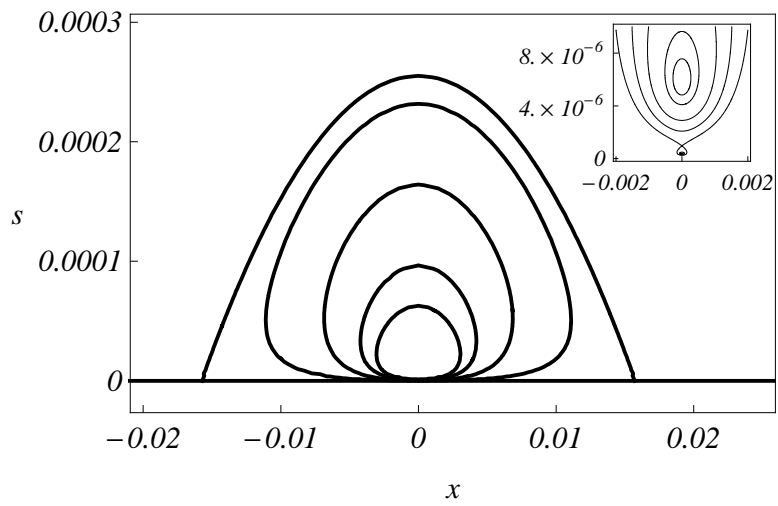}
\caption{Streamlines for Type IV, $\protect\alpha=5,\; \protect\kappa=1.5$,
and $\hat{s}_1\simeq 5.9\times 10^{-6},\,\hat{s}_2\simeq 3.5\times
10^{-7},\, \hat{s} \simeq 9.6\times 10^{-7} ,\, x_{\pm} \simeq \pm 0.015$.}
\label{streamlinecase4}
\end{figure}
\clearpage

\begin{figure}[ht!]
\centering
\includegraphics[width=4in]{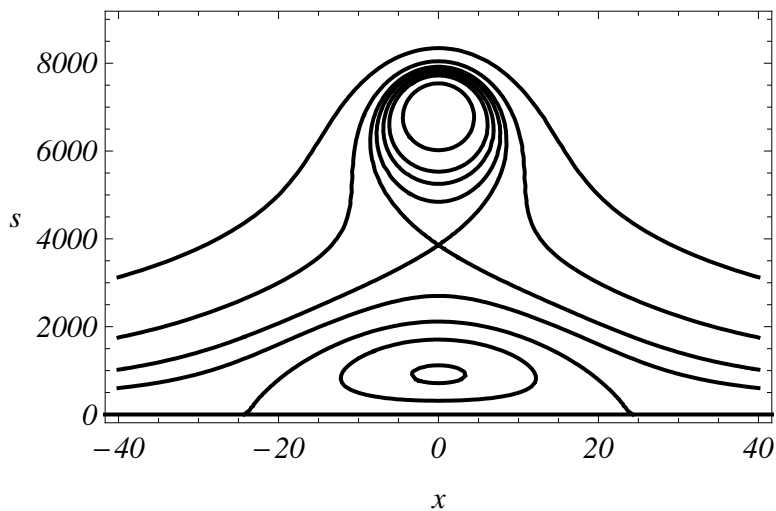}
\caption{Streamlines for Type I, $\protect\alpha=0.1,\; \protect\kappa=1.5$,
and $\hat{s}_1\simeq 895,\,\hat{s}_2\simeq 6914,\, \hat{s} \simeq 3860
,\,x_{\pm} \simeq \pm 24$.}
\label{streamlinehairpin}
\end{figure}

We should point out that for different parameter values of $\alpha ,\,\kappa
$, the structure of the flow could be dramatically different as shown in
Fig. (\ref{streamlinehairpin}) for a Type I configuration for $\alpha $
substantially smaller than unity, where a hairpin vortex structure is
present. We note that in this \textquotedblleft singular\textquotedblright\
case, the dynamics is most certainly inequivalent to any of the instances
illustrated in Figs. 1 - 4. More specifically, the behavior of the stable
and unstable manifolds is very different from those in the other phase
portraits: In this case there is a single homoclinic loop $\ell _{\mathrm{I}%
} $ beginning and ending at $q_{\mathrm{I}}$, which encloses the singular
center $q_{\mathrm{I}_{+}}$. In contrast with the other behaviors
illustrated above, however, the portion of $W^{s}(q_{\mathrm{I}})$ lying
below $q_{\mathrm{II}}$ actually goes off to $x=$ $\infty $, while the
component of $W^{u}(q_{\mathrm{II}})$ below $q_{\mathrm{II}}$ tends to $%
x=-\infty $, with both of these semi-infinite curves remaining above the
cycle $\mathcal{Z}_{\mathrm{I}}$, as depicted in Fig. 5. The analysis of
this singular case - and the bifurcation behavior that transforms the
dynamics in Fig. 5 to that in Fig. 1 - must wait for a later investigation,
as we shall not attempt to treat this in any further detail in the sequel.

We note, however, that hairpin vortex configurations have been observed
computationally in the vortex flow in a cylinder with rotating covers \cite%
{bbs}. When the two covers are co-rotating, a part of the parameter plane
(defined by the Reynolds number and the aspect ratio of the cylinder) may
exist where a hairpin vortex is present in a steady flow. This region of the
parameter plane grows to a significant size when the ratio of the angular
velocities of the covers is increased.

It is eminently clear from our analysis in this section that Type I for a
sufficiently large ambient swirling flow and Type IV , among all four types
treated, produce the type of kinematics most consistent with the B-type
breakdown behavior encountered in experiments and obtained from numerical
solutions of both the Euler and Navier-Stokes equations. Consequently, we
shall only consider Type I with $\alpha \geq 5$ in the sequel (keeping in
mind the close analogs for Type IV) , and take advantage of this by reducing
or dropping all of the unnecessary subscripts and superscripts in order to
simplify our notation. We note only that most of the results that we shall
obtain for this type of coaxial vortex ring positioning - including the
existence of chaotic regimes for slightly oscillating, time-dependent
kinematics - have direct parallels for the other types - and especially Type
IV - as well.

\subsection{The direct approach revisited}

As indicated above, from now on we shall consider Type I (with $\alpha $
large enough to preclude the hairpin vortex behavior shown in Fig. 5) almost
exclusively. So in what follows, we employ simplified notation such as $\hat{%
s}_{1}$ and $\hat{s}_{2}$ for $\hat{s}_{1(\mathrm{I})}$ and $\hat{s}_{2(%
\mathrm{I})}$, respectively, $\hat{s}$ for $\hat{s}_{\mathrm{I}}$, $\nu $
for $\nu _{\mathrm{I}}$, and $x_{\pm }$ for $x_{\mathrm{I}}^{(\pm )}$. Now
we shall look once again at the kinematics from the perspective of the
direct approach formulated in Subsection 3.2, but here for the case when the
rings are located at the stationary points $(\hat{s}_{1},\hat{s}_{2},\xi
,\xi )$ as described above. Since the rings are stationary, the potentially
time-dependent Hamiltonian system (27), or (29)-(30), is actually
autonomous. In particular, owing to our assumption that the $x$-component of
the ambient velocity field is just the constant $-\alpha $, the equations of
motion are
\begin{align}
\dot{s}& =4r\left( x-\xi \right) \sum\limits_{k=1}^{2}\kappa _{k}\hat{r}%
_{k}\int\nolimits_{0}^{\pi /2}\frac{\cos 2\sigma d\sigma }{\left[ \left( r-%
\hat{r}_{k}\right) ^{2}+\left( x-\xi \right) ^{2}+4r\hat{r}_{k}\sin
^{2}\sigma \right] ^{3/2}},  \notag \\
\dot{x}& =-\alpha +2\sum\limits_{k=1}^{2}\kappa _{k}\hat{r}%
_{k}\int\nolimits_{0}^{\pi /2}\frac{\left( \hat{r}_{k}-r\cos 2\sigma \right)
d\sigma }{\left[ \left( r-\hat{r}_{k}\right) ^{2}+\left( x-\xi \right)
^{2}+4r\hat{r}_{k}\sin ^{2}\sigma \right] ^{3/2}}.  \label{e65}
\end{align}

From the dynamical equations (65) in comparison with the last two equations
of (48), we see that the analysis of the stationary points and invariant
manifolds for the passive formulation in this case (in Subsections 5.1 and
5.2) immediately carries over to the direct approach. More specifically, if
the rings are initially at $(\hat{s}_{1},\hat{\xi})$ and $(\hat{s}_{2},\hat{%
\xi})$, respectively, then these points are singular centers, and the other
stationary points are $\left( 0,x_{\pm }\right) $, both of which are saddle
points, and the saddle point $\left( \hat{s},\hat{\xi}\right) $. In
addition, the stable and unstable manifolds of the saddle points are just
the same, with the obvious identifications, as those described above for $%
p_{\pm }$ and $q$, and give rise to the identical heteroclinic cycle - $%
\mathcal{Z}$ - and homoclinic loops $\ell _{+}$ and $\ell _{-}$ forming the
figure eight curve $\zeta $. Accordingly the phase portrait for (27) is
indistinguishable - modulo the obvious change in labels - from that
portrayed in Figs. 1 and 4 for Types I and IV.

It is clear from our analysis so far that for fixed ring configurations, the
dynamics of the passive fluid particles (\emph{i.e}. the kinematics) is
quite regular; in fact, the system is LA-integrable since it is Hamiltonian,
autonomous, and has just one degree-of-freedom. Consequently, the Poincar%
\'{e} map of first returns to the meridian half-plane $\mathfrak{H}$, which
we introduced briefly at the end of Subsection 3.2, cannot admit chaotic
regimes in the fixed ring case, and so must also be regular. Therefore, it
is manifest that if we are to generate chaotic kinematic regimes in our
model - which are of particular interest in this investigation - the coaxial
rings cannot both be stationary. In other words, a necessary condition for
chaotic dynamics of the passive fluid particles is that the ring dynamics be
nontrivial. We shall explore these matters in considerably greater depth in
the next two sections.

\section{Kinematics for Slightly Varying Rings}

A number of studies taken together such as those of Breuer \cite{BREU}, Br\o %
ns \emph{et al}. \cite{BVS, BB, BSSZ}, Gelfgat \emph{et al}. \cite{GBS},
Hartnack et al. \cite{HBS}, Holmes \cite{holmes}, Krause \cite{krause1,
krause2, krause3}, Lopez \& Perry \cite{LP}, Rusak \emph{et al}. \cite{RWW},
Serre \& Bontoux \cite{serbon}, S\o rensen \& Christiansen \cite{SC},
Sotiropoulos \emph{et al}. \cite{sotir1, sotir2}, and Weimer \cite{WEIMER},
present compelling evidence, and in just a few cases proof, that the
kinematics in B-type vortex breakdown exhibits chaotic regimes, which may
include those of Shilnikov type. However, it appears that these chaotic
effects are all traceable to small perturbations that break the axisymmetry
of the configuration. Our model remains perfectly axisymmetric, so
apparently the only way to generate chaotic kinematics is to allow some
non-stationarity for the system. Here we develop a framework for studying
the kinematics of our model with small variations in the ring motion, and
leave to the next section a more detailed analysis of the possible chaotic
regimes generated by slightly varying rings. We begin by analyzing the
dynamics of the two rings for small motions about the fixed points.

\subsection{Dynamics for varying rings}

To represent small motions of the rings about the fixed points, we set%
\begin{equation}
s_{k}=\hat{s}_{k}+\varphi_{k},\;x_{k}=\hat{\xi}+\psi_{k},\quad(1\leq k\leq2)
\label{e66}
\end{equation}
and note that it follows from the LA-integrability of (16), or equivalently
(17)-(18), that the functions $\varphi_{k},\psi_{k}$ in solutions of (16)
are at least quasiperiodic, and are periodic if the initial conditions are
chosen properly. We shall choose the perturbation functions so that $%
\left\vert \varphi_{k}\right\vert $ and $\left\vert \psi_{k}\right\vert $ to
be very small in a manner that shall be specified more precisely in what
follows. The three-dimensional hyperplane
\begin{equation}
G^{-1}\left( \hat{a}\right) =\left\{ \left( s_{1},s_{2},x_{1},x_{2}\right)
:s_{1}+\kappa s_{2}=\hat{s}_{1}+\kappa\hat{s}_{2}:=\hat {a}\right\}
\label{e67}
\end{equation}
of the phase space of (16) is, as pointed out in (26), an invariant manifold
of (16). To simplify matters, we shall only consider (varying) solutions of
(16) that lie in $G^{-1}\left( \hat{a}\right) $, which implies that%
\begin{equation}
\varphi_{2}=-\kappa^{-1}\varphi_{1}.  \label{e68}
\end{equation}
This choice is also motivated by the fact that an examination of (46)
readily reveals that the center surface $W^{c}\left( p\right) $ associated
to any stationary point $p$ of the ring dynamics is contained in $%
G^{-1}\left( \hat{a}\right) $.

If we choose the varying solution (66), with $\left\vert \varphi
_{k}(0)\right\vert ,\left\vert \psi_{k}(0)\right\vert $ $(1\leq k\leq2)$
sufficiently small, so that it is initially on $W^{c}\left( p\right) $, then
our varying rings orbit will be periodic. However, we do not know exactly
what $W^{c}\left( p\right) $ is at this point. Of course, returning once
again to (46) and (47), it is straightforward to show that the linear
approximation $W_{\mathrm{lin}}^{c}\left( p\right) $ of $W^{c}\left(
p\right) $ is defined by the homogeneous linear equations%
\begin{align}
s_{1}+\kappa s_{2} & =0,  \notag \\
Ax_{1}-Bx_{2} & =0,  \label{e69}
\end{align}
where%
\begin{align}
A & =A\left( \hat{s}_{1},\hat{s}_{2},\hat{\xi};\alpha,\kappa\right)
:=\left\{ \nu^{2}+\partial_{x_{1}}\Phi_{2}\left( p\right) \left[
\partial_{s_{2}}\Psi_{1}\left( p\right) -\kappa\partial_{s_{1}}\Psi
_{1}\left( p\right) \right] \right\} ,  \notag \\
B & =B\left( \hat{s}_{1},\hat{s}_{2},\hat{\xi};\alpha,\kappa\right)
:=\partial_{x_{1}}\Phi_{2}\left( p\right) \left[ \partial_{s_{2}}\Psi
_{1}\left( p\right) -\kappa\partial_{s_{1}}\Psi_{1}\left( p\right) \right] ,
\label{e70}
\end{align}
and%
\begin{equation}
\nu=\sqrt{\partial_{x_{1}}\Phi_{2}\left( p\right) \left[ \partial_{s_{2}}%
\Psi_{2}\left( p\right) -\partial_{s_{2}}\Psi_{1}\left( p\right)
+\kappa\left( \partial_{s_{1}}\Psi_{1}\left( p\right)
-\partial_{s_{1}}\Psi_{2}\left( p\right) \right) \right] }.  \label{e71}
\end{equation}
Observe here that we can use the above to approximate the frequency $\omega$
of periodic solutions near the center in $W^{c}\left( p\right) $; namely%
\begin{equation}
\omega\simeq\nu.  \label{e72}
\end{equation}
Moreover, if we select the initial point to be on $W_{\mathrm{lin}%
}^{c}\left( p\right) $ and very near to, but distinct from the stationary
point $p$, we have the following additional useful (but only approximate)
relationship between the perturbing functions - defined in (67) - for the
coordinates $x_{1}$ and $x_{2}$:%
\begin{equation}
\psi_{2}\simeq\frac{A}{B}\psi_{1}.  \label{e73}
\end{equation}

All of the above suggests a simple method for locating a point (distinct
from $p$) of $W^{c}\left( p\right) $: Choose small values for $%
\varphi_{1}\left( 0\right) $ and $\varphi_{2}\left( 0\right) $ that satisfy
(68), and select $\psi_{1}\left( 0\right) $ to be small and nonzero. Then,
if we select a very small values of all of $\varphi_{1}\left( 0\right)
,\varphi_{2}\left( 0\right) $ and $\psi_{1}\left( 0\right) $, we are certain
to find the value of $\psi_{2}\left( 0\right) $ satisfying $\left(
s_{1}(0),s_{2}(0),x_{1}(0),x_{2}(0)\right) \in W^{c}\left( p\right) $ in a
very small interval centered at $\left( A/B\right) \psi_{1}\left( 0\right) $%
; for example, by using a simple bisection method. Naturally, once we have
found such a solution, we are guaranteed that it is periodic and that its
period is approximately equal to $2\pi/\nu$.

For a more precise specification of the size of the perturbations around the
fixed points, we turn to the following result, which is illustrated in Fig.
6.

\bigskip

\noindent \textbf{Lemma 1.} \emph{The varying solutions of }(16) \emph{%
sufficiently close to a point on the center manifold }$W^{c}\left( p\right) $%
\emph{, which is in a sufficiently small neighborhood of the fixed point }$%
p=\left( \hat{s}_{1},\hat{s}_{2},\hat{\xi},\hat{\xi}\right) $\emph{, can be
expressed in terms the following expansions in terms of a small parameter }$%
\mu $:%
\begin{align}
s_{1}& =s_{1}\left( t;\mu \right) =\hat{s}_{1}+\mu \sin \nu t+O\left( \mu
^{2}\right) ,\;s_{2}=\hat{s}_{2}-\left( \mu /\kappa \right) \sin \nu
t+O\left( \mu ^{2}\right) ,  \notag \\
x_{1}& =x_{1}\left( t;\mu \right) =\hat{\xi}+\mu \cos \nu t+O\left( \mu
^{2}\right) ,\;x_{2}=\hat{\xi}+\mu \left( A/B\right) \cos \nu t+O\left( \mu
^{2}\right) ,  \label{e74}
\end{align}%
\emph{where the order relations are uniform in }$t$\emph{\ when }$\left\vert
\mu \right\vert \leq \epsilon _{\ast }:=(1/2)\min \{\hat{s}-\hat{s}_{1},\hat{%
s}_{2}-\hat{s}\}$, \emph{with }$\hat{s}_{1},\hat{s}_{2},$\emph{and }$\hat{s}$
\emph{as defined in our discussion of} (48) \emph{and }(65)\emph{\ in the
preceding section. This restriction guarantees - for one thing - that the
orbits of the rings do not intersect one another.}

\medskip

\noindent\emph{Proof: }This follows directly from the LA-integrability of
the system, the properties of $W_{\mathrm{lin}}^{c}\left( p\right) $ in
relation to $W^{c}\left( p\right) $, and the description of the stationary
points and phase space structure for the corresponding dynamical systems
delineated above. $\blacksquare$


\begin{figure}[th]
\centering
\includegraphics[width=3.3in]{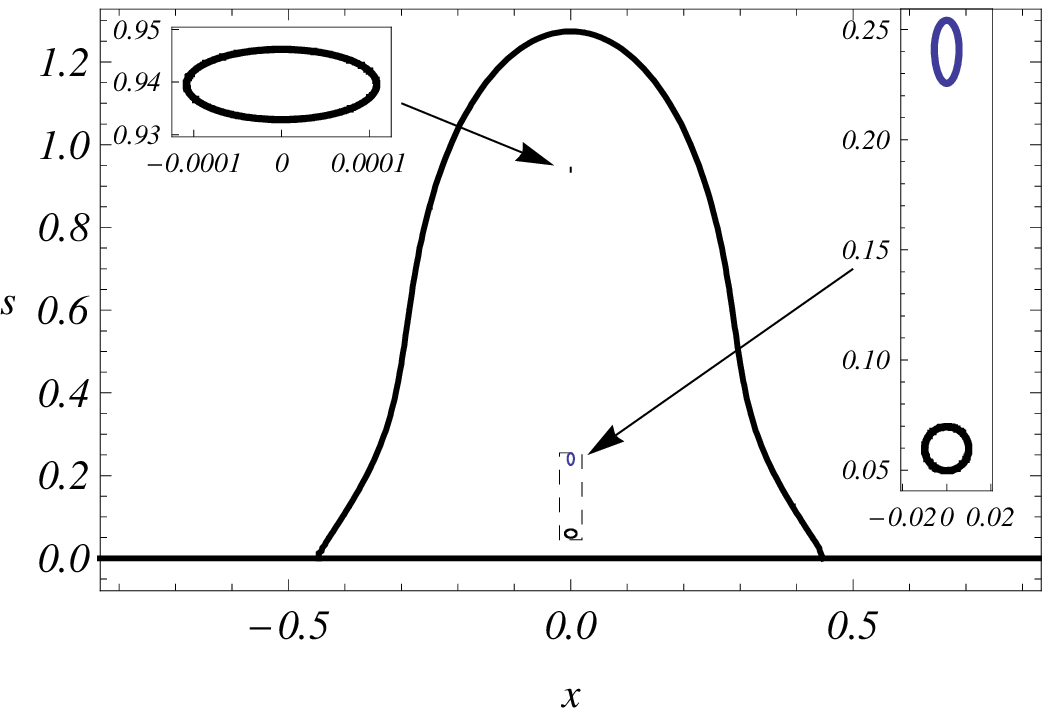}
\includegraphics[width=3.3in]{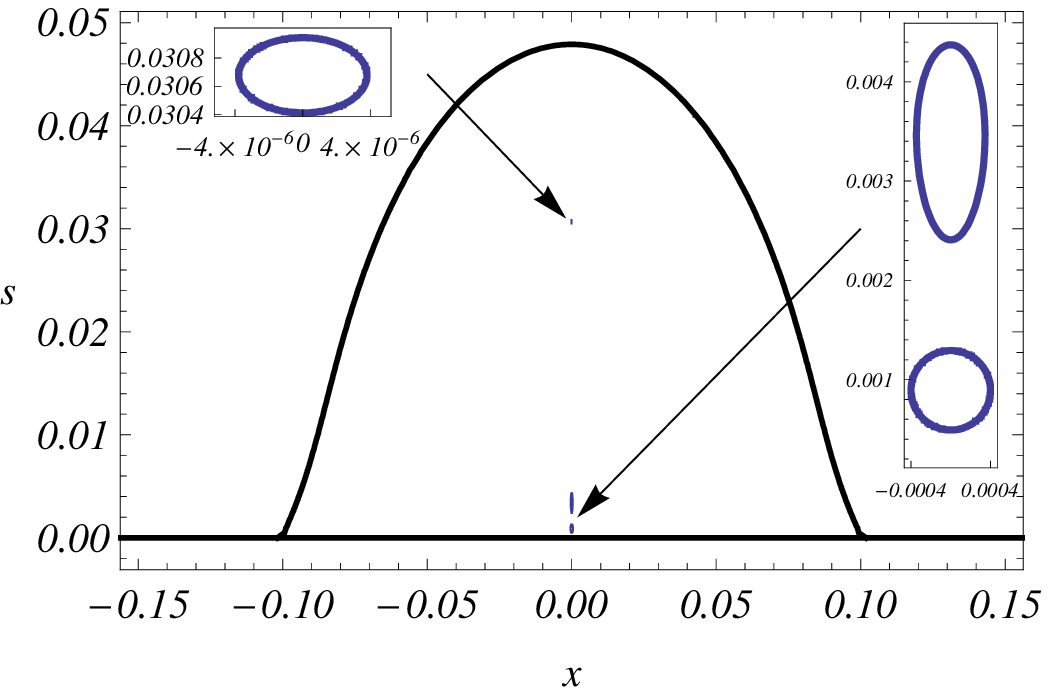}
\caption{Slightly varying ring dynamics and stagnation point inside the
bubble for $\protect\alpha =5,\,\protect\mu =0.01$ (left) and $\protect%
\alpha =20,\,\protect\mu =4\times 10^{-3}$ (right). The unperturbed
heteroclinic orbit is also shown.}
\end{figure}

We note that in Fig. 6 the small rectangular boxes in the upper left hand
corners depict the much smaller scale elliptic orbits of the slightly
oscillating rings in neighborhood of the upper fixed points, while the small
rectangular boxes to the right within the figures provide a microscopic view
of the dynamics near the lower fixed points, which are extremely close to
the axis of symmetry.

\subsection{Kinematics for oscillating rings}

Returning to (29)-(30), it now follows directly from Lemma 1 and a
straightforward calculation, left to the reader, that we can expand (30) in
powers of $\mu$ as described in our next lemma.

\bigskip

\noindent\textbf{Lemma 2.} \emph{For the varying ring dynamics as described
in Lemma 1, we can expand the Hamiltonian function of} (29) \emph{-
representing the direct kinematics formulation - as}%
\begin{equation}
\mathcal{H}\left( s,x,t\right) =\mathcal{H}_{0}\left( s,x\right) +\mu%
\mathcal{H}_{1}\left( s,x,t\right) +O\left( \mu^{2}\right) ,  \label{e75}
\end{equation}
\emph{uniformly in }$t$\emph{\ when }$\left\vert \mu\right\vert \leq
\epsilon_{\ast}$\emph{, where}
\begin{equation}
\mathcal{H}_{0}\left( s,x\right) :=\alpha s-4r\sum\limits_{k=1}^{2}\kappa_{k}%
\hat{r}_{k}\int\nolimits_{0}^{\pi/2}\frac{\cos2\sigma d\sigma}{\hat{\Delta}%
_{k}^{1/2}},  \label{e76}
\end{equation}
\emph{and}%
\begin{align}
\mathcal{H}_{1}\left( s,x,t\right) & :=2r\sin\nu
t\sum\limits_{k=1}^{2}\left( -1\right) ^{k}\int\nolimits_{0}^{\pi/2}\frac{%
\left[ r\left( r-\hat{r}_{k}\right) +\left( x-\hat{\xi}\right) ^{2}+2r\hat{r}%
_{k}\sin ^{2}\sigma\right] \cos2\sigma d\sigma}{\hat{r}_{k}\hat{\Delta}%
_{k}^{3/2}}+  \notag \\
& \qquad\qquad4r\left( x-\hat{\xi}\right) \cos\nu
t\int\nolimits_{0}^{\pi/2}\left\{ \frac{\hat{r}_{1}}{\hat{\Delta}_{1}^{3/2}}+%
\frac{\kappa A\hat{r}_{2}}{B\hat{\Delta}_{2}^{3/2}}\right\} \cos2\sigma
d\sigma.  \label{e77}
\end{align}

\bigskip

\noindent With Lemmas 1 and 2 now at our disposal, we have the elements
necessary to delve more deeply into the behavior of the kinematics - as
characterized by the time-dependent Hamiltonian system (29)-(30).

\subsubsection{Oscillation of stagnation points}

When the rings are stationary, we showed in Subsection 5.1 that the
stagnation points of the advected flow on the $x$-axis are fixed at the
points $x_{\pm}$ determined by (51) and (52). In contrast, when the rings
are oscillating in accord with Lemma 1, these stagnation points will also
vary slightly in a quasiperiodic (or periodic if the rings are properly
tuned) manner. We denote the time-varying positions of the stagnation points
by $x_{\pm}\left( t;\mu\right) $, noting that these are independent of $t$
when $\mu=0$ and given as
\begin{equation}
x_{\pm}\left( t;0\right) =x_{\pm}.  \label{e78}
\end{equation}

These points are determined in a manner analogous to $x_{\pm }$; namely,
they are the solutions in $x$ to%
\begin{equation}
\Psi \left( x,t;\mu ;\kappa \right) =-\alpha +\pi \left\{ \frac{%
r_{1}^{2}\left( t;\mu \right) }{\left[ s_{1}\left( t;\mu \right) +\left(
x-x_{1}(t;\mu )\right) ^{2}\right] ^{3/2}}+\frac{\kappa r_{2}^{2}\left(
t;\mu \right) }{\left[ s_{2}\left( t;\mu \right) +\left( x-x_{2}(t;\mu
)\right) ^{2}\right] ^{3/2}}\right\} =0,  \label{e79}
\end{equation}%
where $s_{k}\left( t;\mu \right) $ and $x_{k}(t;\mu )$, $1\leq k\leq 2$, are
as defined in Lemma 1. It is straightforward to show, just as for the
stationary case described by (51) and (52), that, under the restrictions on $%
\mu $ specified in Lemma 1, (79) has precisely two solutions satisfying
\begin{equation}
x_{-}\left( t;\mu \right) <\hat{\xi}<x_{+}\left( t;\mu \right)  \label{e80}
\end{equation}%
for all $t\in \mathbb{R}$, where $x_{-}\left( t;\mu \right) $ and $%
x_{+}\left( t;\mu \right) $ are only approximately symmetric with respect to
$\hat{\xi}$, rather than exactly so as in the stationary ring case. The
oscillation of the stagnation points is illustrated in Fig. 7 for several
parameter values. Naturally, we would expect this symmetry breaking to have
consequences regarding the nature of the time-varying vortex breakdown
bubble. As we shall see, this loss of axisymmetry can indeed have very
dramatic consequences on the kinematics.


\begin{figure}[ht!]
\centering
\includegraphics[width=3.3in]{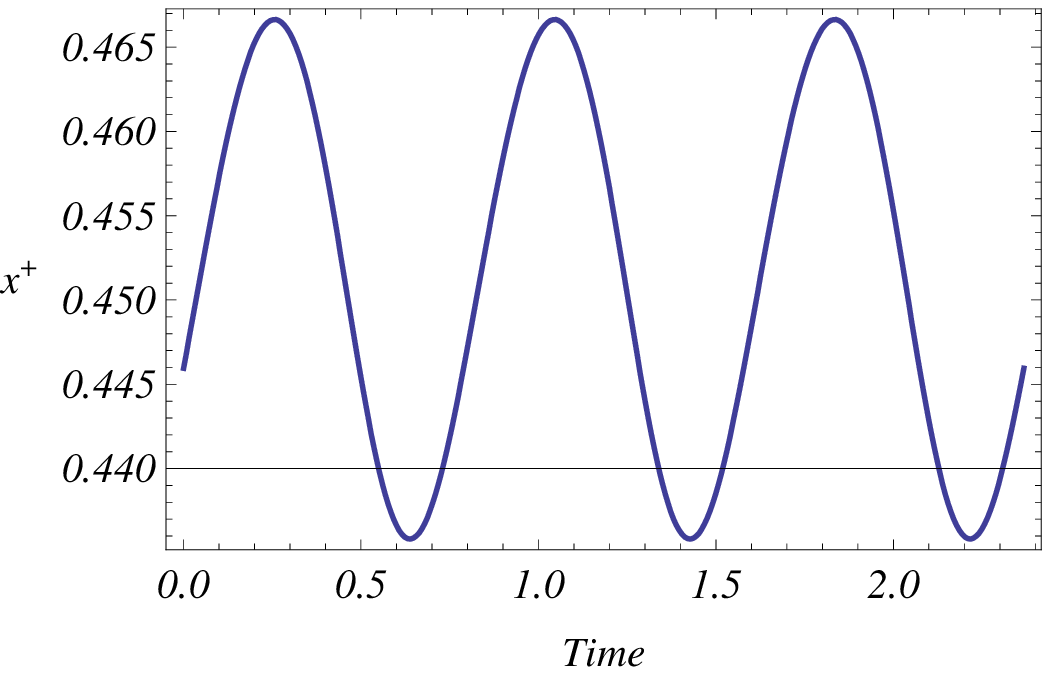}
\includegraphics[width=3.3in]{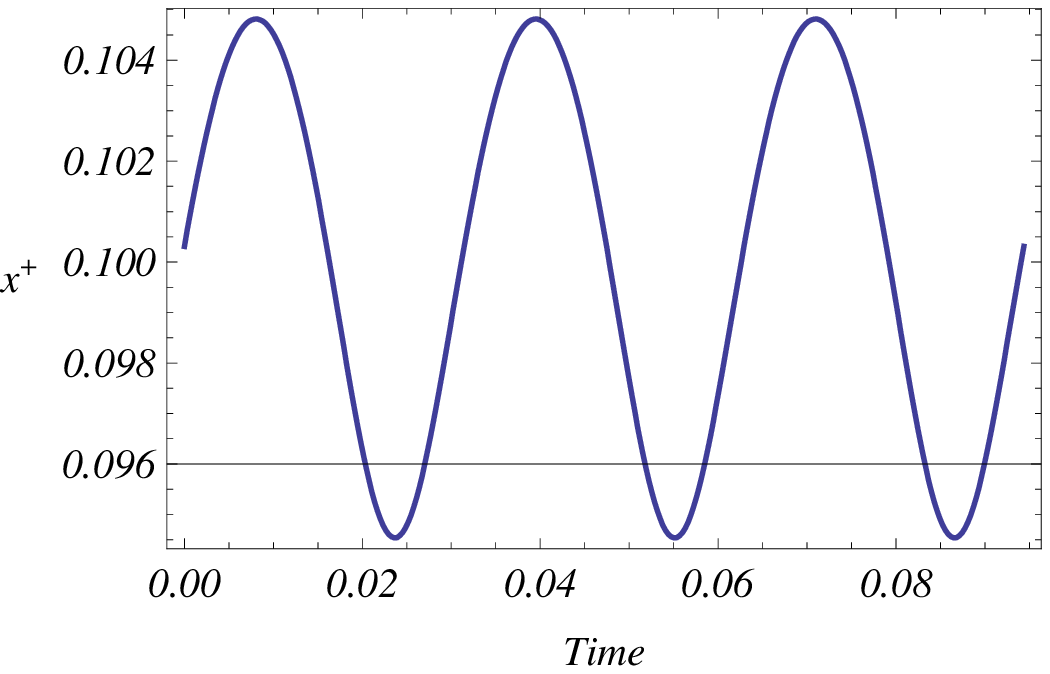}
\caption{Oscillation of the stagnation point along the axis of symmetry ($%
s=0 $) for $\protect\alpha=5,\,\protect\mu=0.01$ (left) and $\protect\alpha%
=20,\, \protect\mu=4\times 10^{-3}$ (right).}
\end{figure}

\subsection{Natural Poincar\'{e} maps for kinematics}

Let us return to the meridional first return Poincar\'{e} map introduced at
the end of Subsection 3.2, which provides a very natural and useful means
for (discretely) characterizing the dynamics of passive fluid particles
induced by the small oscillations in the positions of the coaxial vortex
rings described in (74). For concreteness and simplicity, we shall assume in
the sequel that the rotating component of the superimposed swirling flow has
the (constant) form%
\begin{equation}
\dot{\theta}=\Omega,  \label{e81}
\end{equation}
with $\Omega$ a fixed positive parameter; in other words, $b_{1}=b_{2}=0$ in
(14). In this case, the Poincar\'{e} map has a particularly simple
definition; namely,%
\begin{equation*}
\mathit{\Pi}=\mathit{\Pi}\left( \cdot,\cdot;\alpha,\Omega,\kappa,\mu\right) =%
\mathit{\Pi}_{\mu}:\mathfrak{H}\rightarrow\mathfrak{H},
\end{equation*}
where
\begin{equation}
\mathit{\Pi}_{\mu}\left( s_{0},x_{0}\right) :=\left( s\left( 2\pi
/\Omega\right) ,x\left( 2\pi/\Omega\right) \right) ,  \label{e82}
\end{equation}
and $\left( s\left( t\right) ,x\left( t\right) \right) $ is the unique
solution of (27) satisfying $\left( s\left( 0\right) ,x\left( 0\right)
\right) =\left( s_{0},x_{0}\right) $. Note that although the Poincar\'{e}
map depends on all of the parameters $\alpha,\Omega,\kappa$, and $\mu$ (and
incidentally on $\chi$), we have emphasized its dependence on the parameter $%
\mu$ connected with the oscillation of the vortex rings by singling it out
via a subscript. With this notation, $\mathit{\Pi}_{0}$ corresponds to the
fixed ring case associated to $\mathcal{H}_{0}$ in (75).

As we shall see in our analysis of chaotic motions in the next section using
Melnikov's method, it useful to list a number of elementary properties of
the Poincar\'{e} map $\mathit{\Pi}$, all of which can be readily deduced
from its definition in terms of the flow of the Hamiltonian system
(29)-(30), which is equivalent to (27):

\begin{itemize}
\item[(PM1)] The map $\mathit{\Pi}_{\mu}:\mathfrak{H}\rightarrow\mathfrak{H}$
is a (real) analytic diffeomorphism (actually a symplectomorphism), which is
analytic in both its variables and its parameters when they are restricted
to the range of values delineated above.

\item[(PM2)] $\mathit{\Pi}_{0}$ has hyperbolic fixed points at $\left(
0,x_{+}\right) $ and $\left( 0,x_{-}\right) $ connected by a heteroclinic
cycle $\mathcal{Z}$ comprised of the line segment joining these points and
the curve $\mathcal{C}$ as described in Subsection 5.2. Moreover, it has the
following fixed points: $(\hat{s}_{1},\hat{\xi})$ and $(\hat{s}_{2},\hat{\xi
})$, which are centers of the symplectomorphism $\mathit{\Pi}_{0}$, along
with the saddle point $(\hat{s},\hat{\xi})$ .

\item[(PM3)] If the dynamics of the rings is (nontrivially) periodic of
period $2\pi/\omega$ and $\omega/\Omega=n\in\mathbb{N}$, the set of natural
numbers, both $\left( 0,x_{+}\right) $ and $\left( 0,x_{-}\right) $ are
hyperbolic fixed points of $\mathit{\Pi}_{\mu}$ whenever $\left\vert
\mu\right\vert \leq\epsilon_{\ast}$. On the other hand, if $\Omega/\omega=m$%
, a positive integer greater than one, both $\left( 0,x_{+}\right) $ and $%
\left( 0,x_{-}\right) $ are hyperbolic periodic points of $\mathit{\Pi}%
_{\mu} $ of (least) period $m$.

\item[(PM4)] In any case, if $\mu$ satisfies $\left\vert \mu\right\vert \leq%
\hat{\epsilon}$ for a sufficiently small $0<\hat{\epsilon}\leq
\epsilon_{\ast}$ it follows from the axisymmetry and the standard results on
the persistence of hyperbolic fixed points (see \emph{e.g}. \cite{guho,
wigbook}) that $\mathit{\Pi}_{\mu}$ has a pair of hyperbolic fixed points $%
\tilde{p}_{+}\left( \mu\right) :=\left( 0,\tilde{x}_{+}(\mu)\right) $ and $%
\tilde{p}_{-}\left( \mu\right) :=\left( 0,\tilde{x}_{-}(\mu)\right) $ such
that: (i) $\tilde{x}_{\pm}(0)=x_{\pm}$; (ii) $\tilde{x}_{-}(\mu )<\tilde{x}%
_{+}(\mu)$ for all $\mu$; and (iii) the stable and unstable manifolds of
these fixed points satisfy $W^{s}(\tilde{p}_{+}\left( \mu\right) )=\left\{
\left( 0,x\right) :\tilde{x}_{-}(\mu)<x\right\} $, and $W^{u}(\tilde{p}%
_{-}\left( \mu\right) )=\left\{ \left( 0,x\right) :x<\tilde{x}%
_{+}(\mu)\right\} $. Moreover, $\mathit{\Pi}_{\mu}$ has a hyperbolic fixed
point $\tilde{q}\left( \mu\right) :=\left( \tilde{s}(\mu),\tilde{x}%
(\mu)\right) $ satisfying $\tilde{q}\left( 0\right) :=\left( \hat{s},\hat{\xi%
}\right) $, which remains close to $\left( \hat {s},\hat{\xi}\right) $ for
all $\left\vert \mu\right\vert \leq\hat{\epsilon}$.
\end{itemize}

\noindent Additional insight into the nature of the map $\mathit{\Pi }_{\mu
} $ can be obtained from the results of several numerical experiments
presented in Figs. 8 and 9 in Section 8.

\subsubsection{The heteroclinic cycle}

As we shall see, the heteroclinic cycle associated to the fixed ring
kinematics plays a key role in the identification of chaotic kinematic
regimes for slightly oscillating ring dynamics. Accordingly a more detailed
analysis is in order. To this end, we recall our definition of the
heteroclinic cycle, which joins the leading and trailing points of the
bubble, and encloses all other stationary points along with their stable and
unstable manifolds defining the homoclinic orbits $\ell_{\pm}$ and the
smaller heteroclinic cycle $\zeta$; namely%
\begin{equation}
\mathcal{Z}:=\mathcal{L}_{\pm}\cup\mathcal{C},  \label{e83}
\end{equation}
where $\mathcal{L}_{\pm}$ is the closed axial line segment in $\mathfrak{H}$
defined (for a fixed $\xi=\hat{\xi}$) as $\mathcal{L}_{\pm}:=\left\{ \left(
0,x\right) :x_{-}\leq x\leq x_{+}\right\} $, and $\mathcal{C}$ is the
convex, $\hat{\xi}$-symmetric curve beginning at $\left( 0,x_{+}\right) $
and ending at $\left( 0,x_{-}\right) $ as described in Subsection 5.2.

In virtue of (65) and (PM2), we know that $\mathcal{Z}$ is a heteroclinic
cycle for both the continuous , LA-integrable, Hamiltonian dynamical system
\begin{equation}
\dot{s}=\partial _{x}\mathcal{H}_{0}\left( s,x\right) ,\quad \dot{x}%
=-\partial _{s}\mathcal{H}_{0}\left( s,x\right) ,  \label{e84}
\end{equation}%
where the Hamiltonian function is defined as in (76), and the discrete,
symplectic (integrable) dynamical system defined by the map%
\begin{equation}
\mathit{\Pi }_{0}:\mathfrak{H}\rightarrow \mathfrak{H},  \label{e85}
\end{equation}%
in accordance with (82). To find the heteroclinic orbit in $\mathcal{L}_{\pm
}$, it follows from (65) that we must solve the differential equation%
\begin{equation}
\dot{x}=-\alpha +\pi \left\{ \frac{\hat{r}_{1}^{2}}{\left[ \hat{r}%
_{1}^{2}+\left( x-\hat{\xi}\right) ^{2}\right] ^{3/2}}+\frac{\kappa \hat{r}%
_{2}^{2}}{\left[ \hat{r}_{2}^{2}+\left( x-\hat{\xi}\right) ^{2}\right] ^{3/2}%
}\right\} ,  \label{e86}
\end{equation}%
which is no easy matter. Fortunately, we do not need a complete closed form
solution. To illustrate this, let the heteroclinic orbit satisfying $x(0)=%
\hat{\xi}$ be denoted as $\phi _{l}(t):=\left( s_{l}(t),x_{l}(t)\right) $.
Then, it follows directly from its definition that $s_{l}(t)=0$ for all $%
t\in \mathbb{R}$, and
\begin{equation}
x_{l}(t)=x_{-}+O\left( e^{\varkappa t}\right) \;\mathrm{and}\;\dot{x}%
_{l}(t)=O\left( e^{\varkappa t}\right) \;\mathrm{as}\;t\rightarrow -\infty ,
\label{e87}
\end{equation}%
and%
\begin{equation}
x_{l}(t)=x_{+}+O\left( e^{-\varkappa t}\right) \;\mathrm{and}\;\dot{x}%
_{l}(t)=O\left( e^{-\varkappa t}\right) \;\mathrm{as}\;t\rightarrow \infty ,
\label{e88}
\end{equation}%
where $\varkappa $ is a positive constant.

On the other hand, it is clear from (76) that $\mathcal{C}$ is the solution
curve of the equation%
\begin{equation}
\alpha r-4\sum\limits_{k=1}^{2}\kappa_{k}\hat{r}_{k}\int\nolimits_{0}^{\pi
/2}\frac{\cos2\sigma d\sigma}{\hat{\Delta}_{k}^{1/2}}=0,  \label{e89}
\end{equation}
which is obtained by factoring $r$ out of the energy curve $\mathcal{H}%
_{0}\left( s,x\right) =0$. It is possible through a rather deeper analysis
of (89) to extract many details concerning the form of $\mathcal{C}$, but it
turns out that we just need to use some simple asymptotic representations
for the heteroclinic orbit $\phi_{u}(t):=\left( s_{u}(t),x_{u}(t)\right) $
lying in $\mathcal{C}$ and satisfying the initial condition $%
s_{u}(0)>0,x_{u}(0)=$ $\hat{\xi}$, along with a few elementary symmetry
properties of this orbit. These could be determined by solving the
differential equation obtained from (67) and (89) - a rather formidable
task. Fortunately, the asymptotic and symmetry properties of $\phi_{u}$ that
we require are rather easy to deduce. Owing to the definition of the
heteroclinic orbit, it is easy to verify that%
\begin{equation}
x_{u}(t)=x_{+}+O\left( e^{\varkappa t}\right) \;\mathrm{and}\;\dot{x}%
_{u}(t)=O\left( e^{\varkappa t}\right) \;\mathrm{as}\;t\rightarrow-\infty,
\label{e90}
\end{equation}
and
\begin{equation}
x_{u}(t)=x_{-}+O\left( e^{-\varkappa t}\right) \;\mathrm{and}\;\dot{x}%
_{u}(t)=O\left( e^{-\varkappa t}\right) \;\mathrm{as}\;t\rightarrow\infty,
\label{e91}
\end{equation}
with
\begin{equation}
s_{u}(t)=O\left( e^{\varkappa t}\right) \;\mathrm{and}\;\dot{s}%
_{u}(t)=O\left( e^{\varkappa t}\right) \;\mathrm{as}\;t\rightarrow-\infty,
\label{e92}
\end{equation}
and
\begin{equation}
s_{u}(t)=O\left( e^{-\varkappa t}\right) \;\mathrm{and}\;\dot{s}%
_{u}(t)=O\left( e^{-\varkappa t}\right) \;\mathrm{as}\;t\rightarrow\infty,
\label{e93}
\end{equation}
where the positive constant $\varkappa$ may have to be adjusted in order to
satisfy all six expressions (88)-(89) and (90)-(93). As for the symmetry
properties, it is clear from the relevant definitions that the following
obtain:%
\begin{equation}
s_{u}(-t)=s_{u}(t)>0\;\mathrm{and}\;x_{u}(-t)-\hat{\xi}=-(x_{u}(t)-\hat{\xi }%
),(\mathrm{with\;}(x_{u}(t)-\hat{\xi})<0\;\mathrm{for}\;t>0)  \label{e94}
\end{equation}
and%
\begin{equation}
\dot{s}_{u}(-t)=-\dot{s}_{u}(t),(\mathrm{with}\;\dot{s}_{u}(t)<0\;\mathrm{%
for\;}t>0)\;\mathrm{and}\;\dot{x}_{u}(-t)=\dot{x}_{u}(t)<0  \label{e95}
\end{equation}
for all $-\infty<t<\infty$.

\section{Chaotic Kinematics}

We shall demonstrate in this section that small oscillations in the coaxial
vortex ring motions lead to chaotic dynamical regimes for the advected fluid
particles: in short, small oscillations generate chaotic kinematics. This
shall be accomplished by using Melnikov's method to demonstrate that slight
motions of the rings can lead to transverse intersections in the
heteroclinic cycle $\mathcal{Z}$, which generate chaotic regimes in virtue
of well known results such as presented in Wiggins \cite{wigbook}. We note
that the axisymmetry of the dynamics and kinematics precludes such
intersections in the $\mathcal{L}_{\pm}$ portion of the cycle, so the only
possible such intersections must occur between the unstable manifold $W^{u}(%
\tilde{p}_{+}\left( \mu\right) )$ and stable manifold $W^{s}(\tilde{p}%
_{-}\left( \mu\right) )$, which comprise $\mathcal{C}$ in the stationary
ring configuration for $\mu=0$.

The chaotic behavior can be summarized in the following result, which also
can be proved by making a fairly straightforward - but far from obvious -
modification of the three coaxial ring analysis in Bagrets \& Bagrets \cite%
{bb}.

\bigskip

\noindent\textbf{Theorem 1}.\emph{\ For each }$0<\left\vert \mu\right\vert
\leq\hat{\epsilon}$\emph{, the map}
\begin{equation*}
\mathit{\Pi}_{\mu}:\mathfrak{H}\rightarrow\mathfrak{H},
\end{equation*}
\emph{defined by }(82)\emph{\ has a heteroclinic cycle }$\mathcal{Z}_{\mu}$%
\emph{\ (with }$\mathcal{Z}_{0}=\mathcal{Z}$\emph{) having transverse
intersections, which implies the existence of chaotic orbits.}

\emph{\medskip}

\noindent\emph{Proof}: The key to the proof is an analysis of the zeros of
the Melnikov function (\emph{cf}. \cite{bb}, \cite{dbiutam}, \cite{blkn2},
\cite{kathas}, \cite{malwig}, \cite{MSWI}, \cite{wigbook}, and \cite{ziglin}%
)
\begin{equation}
\mathfrak{M}\left( \tau\right) =\mathfrak{M}\left( \tau;\alpha
,\Omega,\kappa,\chi\right) :=\int\nolimits_{-\infty}^{\infty}\left\{
\mathcal{H}_{0},\mathcal{H}_{1}\right\} _{0}\left( \phi_{u}(t),t+\tau
\right) dt,  \label{e96}
\end{equation}
where the standard Poisson bracket is defined as usual as
\begin{equation}
\left\{ \mathcal{H}_{0},\mathcal{H}_{1}\right\} _{0}:=\partial _{x}\mathcal{H%
}_{0}\partial_{s}\mathcal{H}_{1}-\partial_{s}\mathcal{H}_{0}\partial_{x}%
\mathcal{H}_{1}.  \label{e97}
\end{equation}
It is clear from (29) that the Melnikov function may be written as%
\begin{equation}
\mathfrak{M}\left( \tau\right) =\int\nolimits_{-\infty}^{\infty}\left[ \dot{s%
}_{u}(t)\partial_{s}\mathcal{H}_{1}\left( \phi_{u}(t),t+\tau\right) +\dot{x}%
_{u}(t)\partial_{x}\mathcal{H}_{1}\left( \phi_{u}(t),t+\tau\right) \right]
dt,  \label{e98}
\end{equation}
and this can readily be shown to be well defined owing to the convergent
nature of the integral, which follows directly from (84) and (90)-(93).

Now we compute from (77) that%
\begin{equation}
\partial_{s}\mathcal{H}_{1}\left( \phi_{u}(t),t+\tau\right) =\Theta
_{1}\left( t\right) \left( \sin\nu t\cos\nu\tau+\cos\nu t\sin\nu \tau\right)
+\Theta_{2}\left( t\right) \left( \cos\nu t\cos\nu\tau-\sin\nu
t\sin\nu\tau\right) ,  \label{e99}
\end{equation}
and%
\begin{equation}
\partial_{x}\mathcal{H}_{1}\left( \phi_{u}(t),t+\tau\right) =2r_{u}\left(
t\right) \left[ \Theta_{3}\left( t\right) \left( \sin\nu t\cos\nu
\tau+\cos\nu t\sin\nu\tau\right) +\Theta_{4}\left( t\right) \left( \cos\nu
t\cos\nu\tau-\sin\nu t\sin\nu\tau\right) \right] ,  \label{e100}
\end{equation}
where%
\begin{align}
\Theta_{1}\left( t\right) & :=\sum\limits_{k=1}^{2}\left( -1\right) ^{k}\hat{%
r}_{k}^{-1}\int\nolimits_{0}^{\pi/2}\hat{\Delta}_{k}^{-5/2}\left(
\phi_{u}(t)\right) \left\{ r_{u}^{-1}(t)\left( x_{u}(t)-\hat{\xi}\right) ^{2}%
\left[ \left( r_{u}(t)-\hat{r}_{k}\right) ^{2}+\left( x_{u}(t)-\hat{\xi}%
\right) ^{2}+\right. \right.  \notag \\
& \qquad\left. \left. r_{u}(t)\hat{r}_{k}\left( 2-\cos2\sigma\right) \right]
-2\hat{r}_{k}\cos2\sigma\left[ \left( r_{u}(t)-\hat{r}_{k}\right)
^{2}-2r_{u}(t)\hat{r}_{k}\sin^{2}\sigma\right] \right\} \cos2\sigma d\sigma,
\label{e101}
\end{align}%
\begin{align}
\Theta_{2}\left( t\right) & :=2\left( x_{u}(t)-\hat{\xi}\right)
\int\nolimits_{0}^{\pi/2}\left\{ \hat{r}_{1}r_{u}^{-1}(t)\hat{\Delta}%
_{1}^{-5/2}\left( \phi_{u}(t)\right) \left[ \left( x_{u}(t)-\hat{\xi }%
\right) ^{2}-\left( r_{u}(t)-\hat{r}_{1}\right) \left( 3r_{u}(t)+\hat {r}%
_{1}\right) -\right. \right.  \notag \\
& \quad\qquad\left. \left. 2r_{u}(t)\hat{r}_{1}\sin^{2}\sigma\right] +\frac{%
\kappa A}{B}\hat{r}_{2}r_{u}^{-1}(t)\hat{\Delta}_{2}^{-5/2}\left(
\phi_{u}(t)\right) \left[ \left( x_{u}(t)-\hat{\xi}\right) ^{2}-\right.
\right.  \notag \\
& \qquad\qquad\qquad\left. \left. \left( r_{u}(t)-\hat{r}_{2}\right) \left(
3r_{u}(t)+\hat{r}_{2}\right) -2r_{u}(t)\hat{r}_{2}\sin^{2}\sigma\right]
\right\} \cos2\sigma d\sigma,  \label{e102}
\end{align}%
\begin{align}
\Theta_{3}\left( t\right) & :=\left( x_{u}(t)-\hat{\xi}\right)
\sum\limits_{k=1}^{2}\left( -1\right) ^{k}\hat{r}_{k}^{-1}\int
\nolimits_{0}^{\pi/2}\hat{\Delta}_{k}^{-5/2}\left( \phi_{u}(t)\right) \left[
2r_{u}(t)\hat{r}_{k}\sin^{2}\sigma-\right.  \notag \\
& \quad\quad\left. \left( r_{u}(t)-\hat{r}_{k}\right) \left( r_{u}(t)+2\hat{r%
}_{k}\right) -\left( x_{u}(t)-\hat{\xi}\right) ^{2}\right] \cos2\sigma
d\sigma,  \label{e103}
\end{align}
and%
\begin{align}
\Theta_{4}\left( t\right) & :=2\int\nolimits_{0}^{\pi/2}\left\{ \hat {r}_{1}%
\hat{\Delta}_{1}^{-5/2}\left( \phi_{u}(t)\right) \left[ \left( r_{u}(t)-\hat{%
r}_{1}\right) ^{2}-2\left( x_{u}(t)-\hat{\xi}\right) ^{2}+4r_{u}(t)\hat{r}%
_{1}\sin^{2}\sigma\right] +\right.  \notag \\
& \quad\left. \frac{\kappa A}{B}\hat{r}_{2}\hat{\Delta}_{2}^{-5/2}\left(
\phi_{u}(t)\right) \left[ \left( r_{u}(t)-\hat{r}_{2}\right) ^{2}-2\left(
x_{u}(t)-\hat{\xi}\right) ^{2}+4r_{u}(t)\hat{r}_{2}\sin^{2}\sigma\right]
\right\} \cos2\sigma d\sigma.  \label{e104}
\end{align}
Observe that it follows directly from (94) and (95) that $\Theta_{1}$ and $%
\Theta_{3}$ are even functions of $t$ while $\Theta_{2}$ and $\Theta_{4}$
are odd functions.

Our analysis so far, together with the identification of odd and even
functions, leads immediately to the following simplification of the Melnikov
function:%
\begin{align}
\mathfrak{M}\left( \tau\right) & =2\cos\nu\tau\int\nolimits_{0}^{\infty
}\left\{ \left[ \dot{s}_{u}(t)\Theta_{1}(t)+2r_{u}(t)\dot{x}%
_{u}(t)\Theta_{3}\left( t\right) \right] \sin\nu t+\right.  \notag \\
& \qquad\quad\left. \left[ \dot{s}_{u}(t)\Theta_{2}(t)+2r_{u}(t)\dot{x}%
_{u}(t)\Theta_{4}\left( t\right) \right] \cos\nu t\right\} dt.  \label{e105}
\end{align}
Whereupon, a straightforward but rather tedious analysis of the integral
shows that
\begin{equation}
\mathfrak{M}\left( \tau\right) =C\cos\nu\tau,  \label{e106}
\end{equation}
where $C>0$. Consequently, $\mathfrak{M}\left( \tau\right) $ has simple
zeros, which proves that there is a transverse heteroclinic orbit in the
cycle $\mathcal{Z}_{\mu}$, and this implies the existence of chaotic
dynamics for the discrete system generated by $\mathit{\Pi}_{\mu}$ for $%
\mu\neq0$. Thus the proof is complete. $\blacksquare$

\bigskip

A few remarks are in order concerning the nature of the chaotic kinematics
described in Theorem 1. We have shown that small oscillations of the rings
leads to a transverse intersection in the large heteroclinic cycle $\mathcal{%
Z}$, but one can also prove that such motions induce transverse
intersections in the homoclinic loops $\ell _{\pm }$ comprising the figure
eight curve $\zeta $(cf. \cite{bb}). Consequently, it is not surprising that
there is an accumulation of chaotic streamlines around the stationary saddle
point $(\hat{s},\hat{\xi})$ as seen in several of our figures. Of course,
the chaotic behavior we have proved cannot - in virtue of axisymmetry -
include a regime of the type discovered by Shilnikov \cite{shil1} and which
has now become a standard fixture in modern dynamical systems theory (\emph{%
cf}. \cite{guho, kathas, wigbook}). However, a close look at the chaos
resulting from a transverse intersection of the stable and unstable
manifolds comprising the upper branch $\mathcal{C}$ of the heteroclinic
cycle $\mathcal{Z}$, reveals certain Shilnikov-like features of the dynamics
of the passive fluid particles. In particular, denumerably many orbits
through points arbitrarily close to the axis of symmetry $\mathcal{L}_{\pm }$
(which contains the unstable manifold of the trailing point $\tilde{x}%
_{-}(\mu )$ described in (PM1)-(PM4)) of the bubble, accumulate chaotically
around the trailing and leading points of the bubble, which is a
characteristic and fundamental property of the dynamics of the
three-dimensional Shilnikov model. We note that dynamical systems
considerations and several experimental and numerical studies of fully
three-dimensional B-type vortex breakdown regimes subject to
non-axisymmetric perturbations indicate that Shilnikov chaos is quite common
in the streamline patterns within the bubbles \cite{BSSZ, holmes, sotir1,
sotir2}.

\section{Illustrative Examples}

In this section we study several cases that demonstrate the kinematic
behavior of our model (for the Type I positioning of the rings) via the
Poincar\'{e} map $\mathit{\Pi }_{\mu }$, and provide insights into the ways
in which it depends on the various parameters, especially with regard to the
onset of chaotic regimes. We shall for each of the cases take $\chi =1,000$,
which is tantamount to prescribing the core radius of the rings to be $%
O(0.001)$. To give a fairly representative range of possibilities, we shall
in particular compute and present pictorial representations of the dynamics
of $\mathit{\Pi }_{\mu }$ for the following cases:

\bigskip


\textbf{Case 1:} $\alpha =5,\Omega =\nu \simeq 7.96,\kappa =1.5$: \textbf{%
(a): }$\mu =0$; \textbf{(b):} $\mu =0.001$; \textbf{(c):} $\mu =0.01$.

\bigskip

\textbf{Case 2:} $\alpha =20,\Omega =\nu \simeq 199.7,\kappa =1.5$: \textbf{%
(a):} $\mu =0$; \textbf{(b):} $\mu =4\times 10^{-5}$; \textbf{(c):} $\mu
=4\times 10^{-4}$

\newpage

\begin{figure}[ht!]
\centering
\includegraphics[width=3.3in]{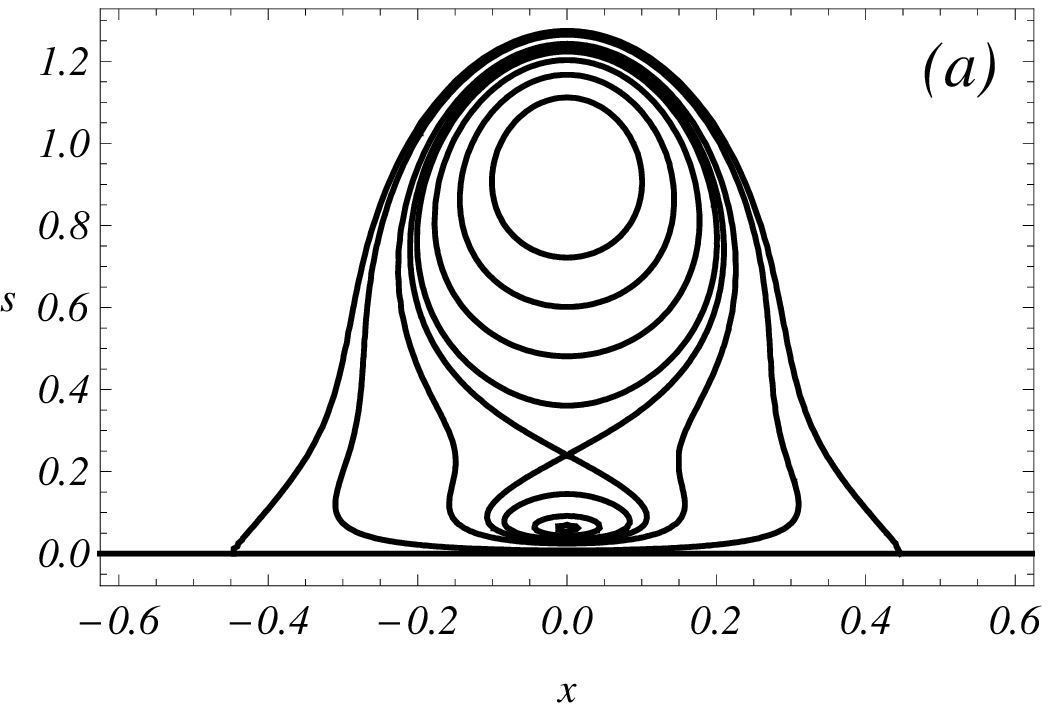}%
\newline
\includegraphics[width=3.3in]{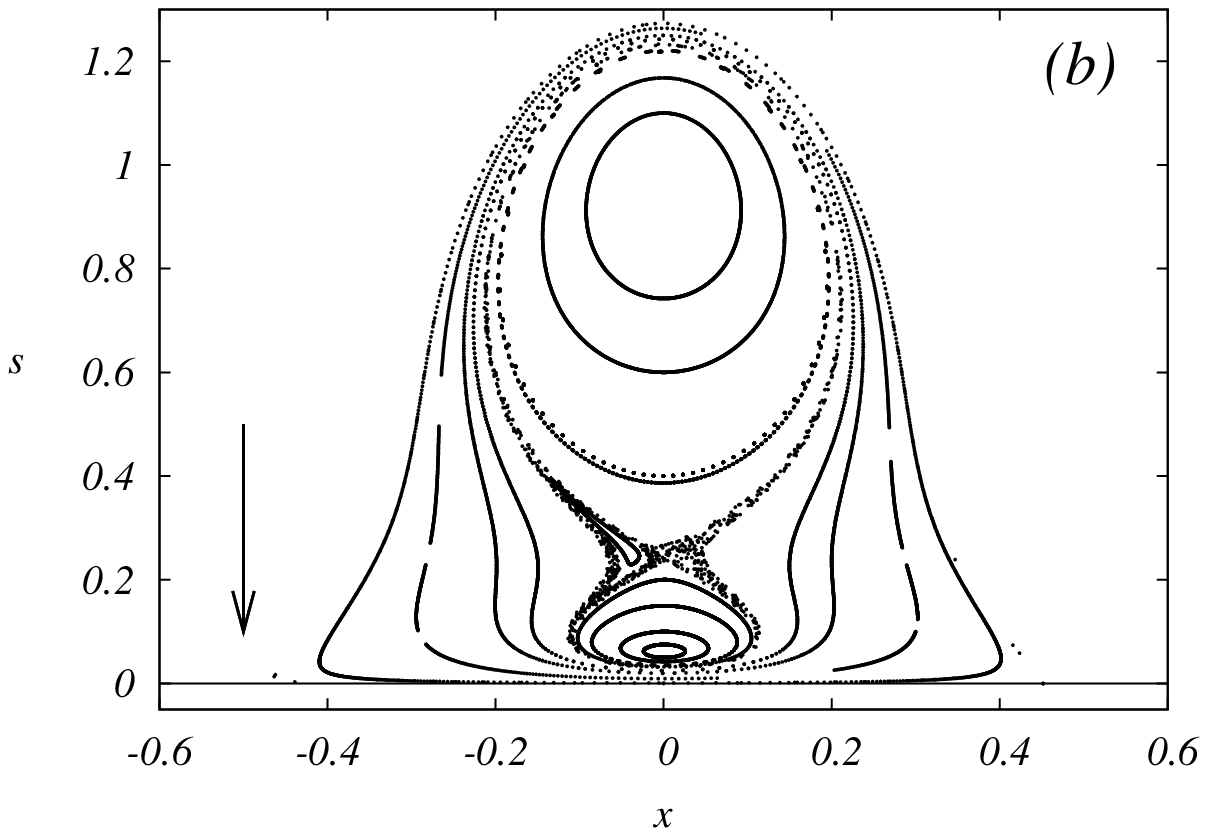} %
\includegraphics[width=3.3in]{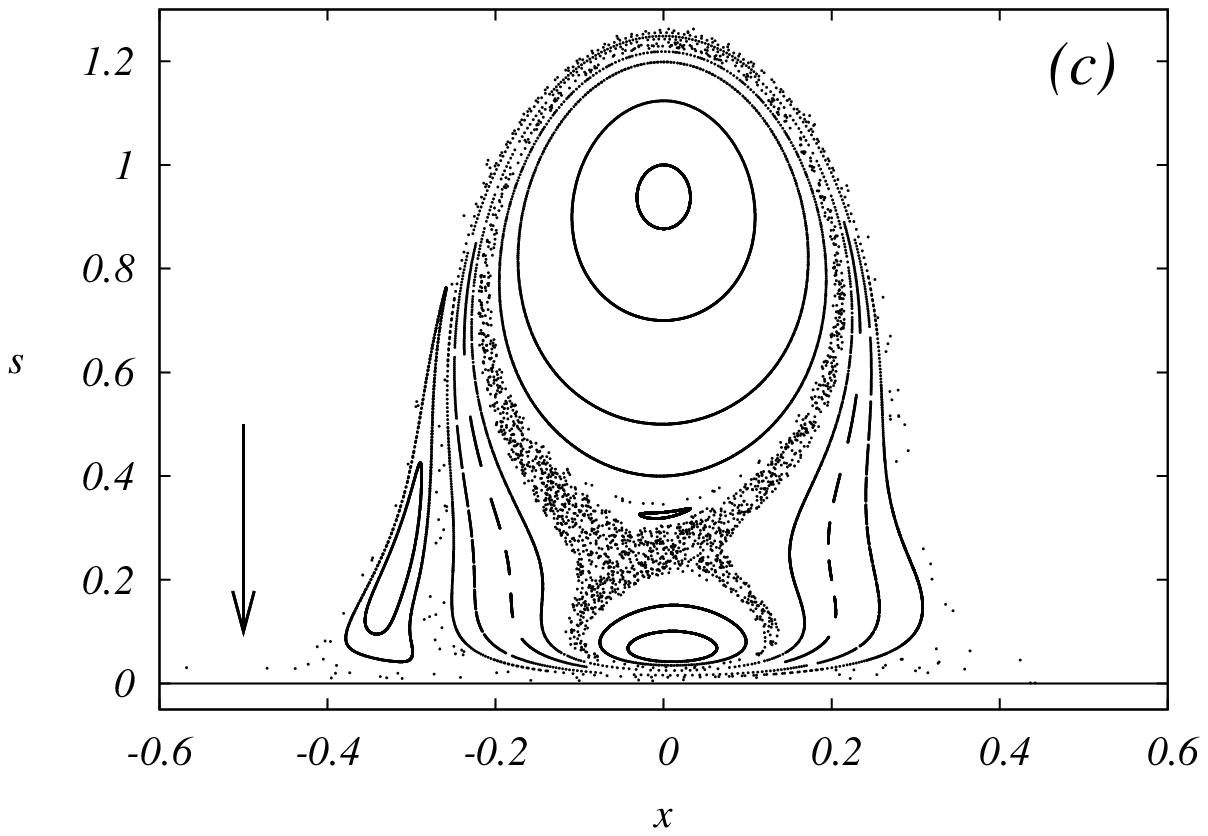}
\caption{Poincar\'e section for $\protect\alpha=5,\, \protect\kappa=1.5,
\Omega=\protect\nu\simeq 7.96$ for a) $\protect\mu=0$, b) $\protect\mu=0.001$
and c) $\protect\mu=0.01$.}
\label{poincarealpha5}
\end{figure}

\newpage

\begin{figure}[ht!]
\centering
\includegraphics[width=3.3in]{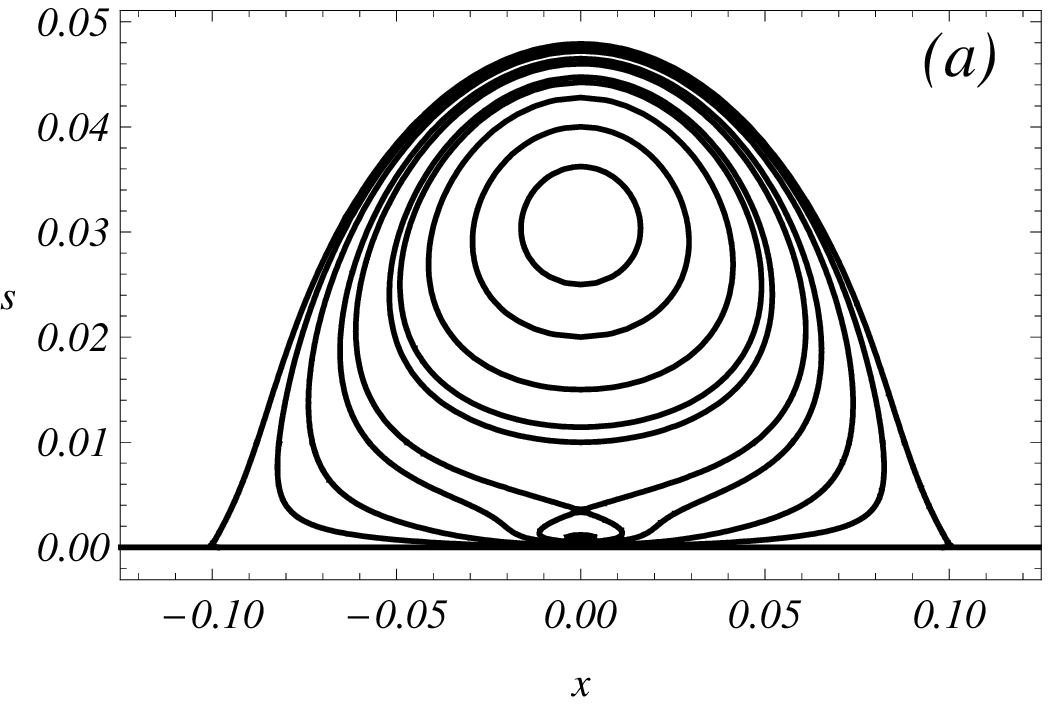}%
\newline
\includegraphics[width=3.3in]{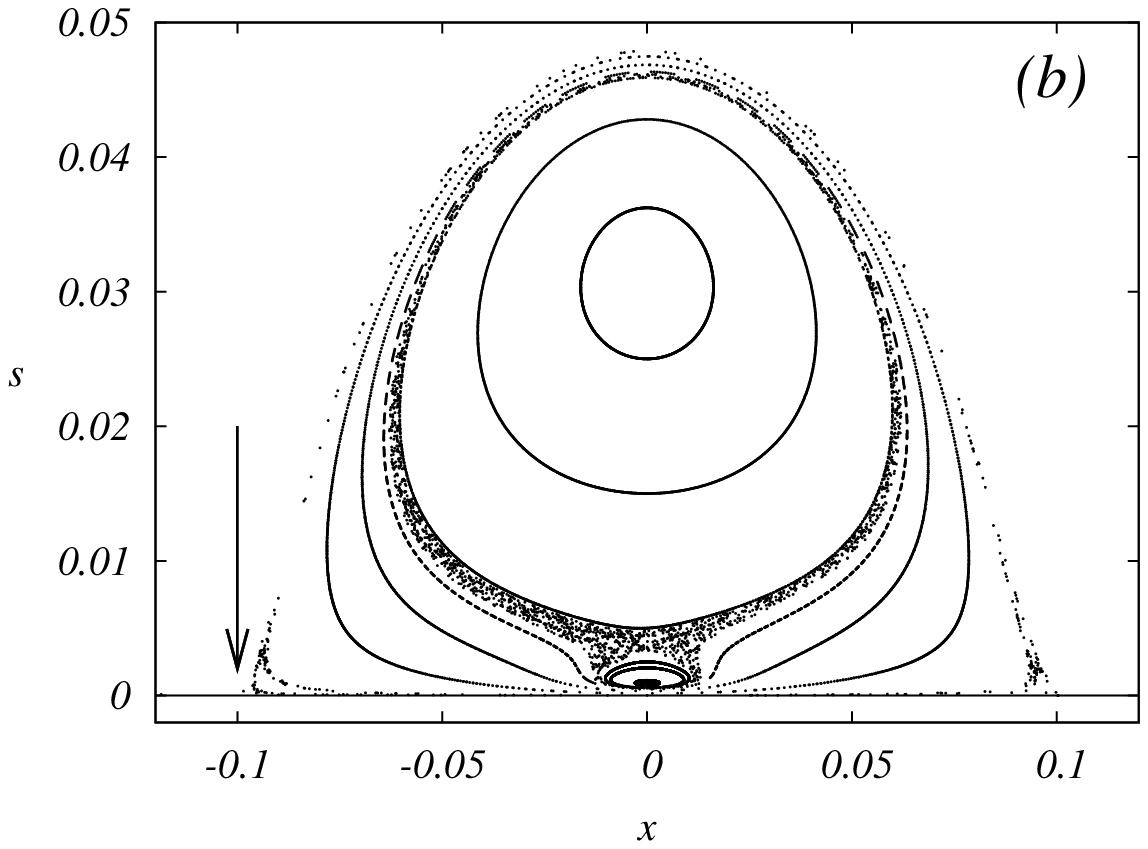}
\includegraphics[width=3.3in]{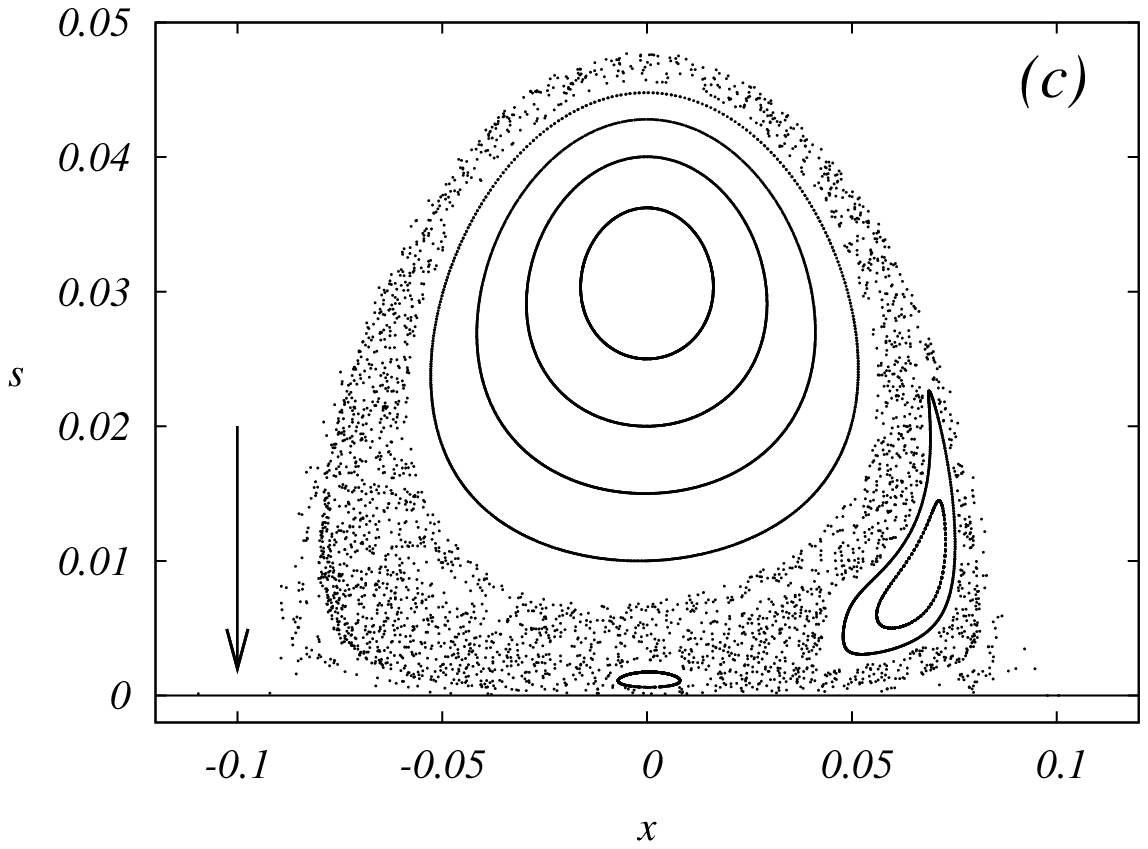}
\caption{Poincar\'e section for $\protect\alpha=20,\, \protect\kappa=1.5,
\Omega=\protect\nu\simeq 199.7$ for a) $\protect\mu=0$, b) $\protect\mu%
=4\times 10^{-5}$ and c) $\protect\mu=4\times 10^{-4}$.}
\label{poincarealpha20}
\end{figure}

The Hamilton equations are solved using a 4th order Runge-Kutta scheme with
a time step of the order of $10^{-5}$.\newline
From the result of the previous section, one expect the homoclinic and
heteroclinic orbit to break after a small perturbation is added to the
system as seen in Figures \ref{poincarealpha5} and \ref{poincarealpha20}. We
see from these figures that as the size of the perturbation increases,
transverse homoclinic intersections interior to the bubble and transverse
heteroclinic cycle behavior of the bubble boundary becomes more pronounced,
so more regular tori are destroyed, leading to larger chaotic regions inside
the bubble - indicated by the characteristic splattering, along with islands
in a sea of chaos. Note in particular how the heteroclinic orbit delineating
the upper boundary of the bubble appears to create a chaotic layer
surrounding the bubble, and that the thickness of this layer increases with
the size of the amplitude of the ring oscillations. One also should notice
the escape of passive particles at the tails of the bubbles indicated by the
arrows, which presumably occurs in devil's staircase fashion (\emph{cf}.
\cite{BSSZ, HBS, sotir2})

\newpage

\section{Epilogue}

We have provided a rigorous demonstration that a pair of coaxial vortex
rings immersed in a swirling ideal fluid flow appear to generate all of the
kinematic behavior associated with vortex breakdown of B-type - including
chaotic streamlines with some Shilnikov-like features when the rings undergo
small oscillations. Naturally, although our results may turn out to be very
useful in studying vortex breakdown, we freely admit that we have by no
means proved that multiple coaxial vortex rings are the fundamental
underlying mechanism in vortex breakdown phenomena. But even if the multiple
rings do turn out to be the ultimate engine of vortex breakdown, one may
well ask - as we have asked ourselves - how are these rings formed in
swirling flows? The answer to this question might turn out to be quite
elusive; nevertheless, we shall next endeavor to posit a plausible
explanation that combines some mathematics and basic fluid mechanics.

To see how the two slender coaxial vortex ring structure (described in Types
I - IV, and perhaps most aptly for Type IV configurations) might develop, we
start with only the swirling flow having its vorticity (vector) along the
axis of symmetry, which we have taken to be the $x$-axis in our analysis. At
the instant when a breakdown bubble is formed, the projected flow on any
meridian half-plane must have at least one stationary point within the
bubble owing to the Poincar\'{e}-Hopf index theorem \cite{arnh, guho,
kathas, wigbook}. As the flow structure of the bubble persists even in
non-Hamiltonian real (viscous and compressible) fluid dynamics, this
observability of the B-type breakdown phenomenon is suggestive of structural
stability for the bubble formation process, which requires that the initial
stationary point must be a saddle, associated, in virtue of the axisymmetry,
with the intersection of any meridian half-plane with a circular streamline
about the $x$-axis corresponding to a hyperbolic periodic orbit of the flow.
But then, once again by the Poincar\'{e}-Hopf index theorem, at least two
centers are required to balance the saddle point. Inasmuch as the process of
bubble formation is essentially instantaneous, one would expect the
necessary two centers to be singular, which is consistent with the
appearance of a pair of coaxial vortex rings. There is also a vorticity
conservation element in this process in accord with Kelvin's theorem \cite%
{chormars, lamb, new, saff, lurupbook}, which may be explained along the
following lines: The appearance of the bubble alters the velocity field in a
way that introduces a vorticity field transverse to the initial field along
the $x$-axis. The induced transverse vorticity component that results from
the bubble formation must be balanced by a structure within the bubble such
as a coaxial vortex ring of nonzero strength, and as we have seen there must
be two rings in order to satisfy the index theorem.

Some of the flow structures found here are richer than those found
computationally for the flow in a cylinder with rotating covers \cite%
{bbs,BVS}. The general appearance of a hyperbolic closed streamline found in
the present study gives rise to structures which, borrowing a notion from
\cite{BVS}, are \textquotedblleft bubbles with an inner
structure\textquotedblright\ as shown in Figs.~1-4. These are only rarely
found in the flows in confined cylinders. However, as the hyperbolic closed
streamline occurs quite close to the axis, the overall flow patterns do not
deviate much from those observed in real flows, even if from a strict
topological perspective they are quite different.

Although our scenario for two coaxial vortex rings driving the B-type vortex
breakdown is reasonably consistent with the governing physical laws and
observed behaviors, we have no illusions about it being indisputable. Be
that as it may, we have shown that such a model does provide a rather
complete description of B-type breakdown flows. Accordingly our approach has
the potential for providing an excellent paradigm and vehicle for developing
techniques and instrumentalities for ameliorating, focusing and controlling
B-type vortex breakdown phenomena, and we hope to investigate some possible
applications in these directions in the not too distant future.

During the course of the analysis of our two ring model, we were struck by
the strong similarities between the kinematics induced in a meridian
half-plane and the phase plane structure of the kinematics generated by a
pair of point vortices in a half-plane. This is hardly surprising in view of
the strong connections between coaxial ring dynamics and point vortex
dynamics established in such investigations as \cite{blkn1}, \cite{BWC}, and
\cite{btk}. Moreover, it has been proved \cite{dbiutam} that the motion of a
pair of point vortices in a half-plane is capable of generating chaotic
behavior in the advection of passive fluid particles. Given the simplicity
of the kinematics generated by a pair of point vortices in a half-plane
compared to that of a pair of coaxial vortex rings, which is essentially on
the order of logarithmic versus elliptic integral Hamiltonians, it seems
natural to investigate the possibility of the former - perhaps immersed in a
parallel flow - producing the complete spectrum of B-type vortex breakdown
kinematics. The practicability of employing the two point vortex model as a
more tractable paradigm for vortex breakdown phenomena of B-type is quite an
attractive option, which definitely calls for further exploration.

Speaking of natural problems and questions for further study, three that we
have avoided in this paper come to mind. First, ambient swirling flows with
translational and rotational components that are dependent on the radial
distance from the axis of symmetry are certainly more representative of the
behavior of real fluids, where one expects to see velocity and rotation
profiles that tend to decrease with distance from the axis. Models
incorporating such radial dependence should be investigated, and one expects
to be able to obtain even more realistic B-type breakdown behavior in such
cases. Secondly, although we pointed out the possibility of rather
unexpected hairpin vortex + bubble structures, seeming to appear from
\textquotedblleft out of the blue\textquotedblright , associated with Type I
for small translational velocities in the ambient swirling flow, we made no
attempt in this paper to analyze them. There is clearly some very
interesting - and challenging - bifurcation behavior associated with these
unusual structures warranting further investigation. In particular, how does
the hairpin + bubble geometry bifurcate into a single bubble as the ambient
translational velocity increases? Lastly, it might be interesting to study
analogous vortex breakdown models for bubble structures based upon three or
more slender coaxial vortex rings immersed in swirling flows, again in aid
of obtaining even more realistic B-type vortex breakdown structures.

Perfect axisymmetry is, of course, as much an abstract ideal as is an ideal
fluid, and as we have already remarked, there have been a variety of
non-axisymmetric behaviors observed in experiments and numerical studies,
which are precluded by the enforced axisymmetry of our model. It is
therefore quite natural to enquire whether or not a full three-dimensional
version of our coaxial vortex ring model that admits small symmetry breaking
perturbations is capable of producing all of these non-axisymmetric
kinematic effects? We plan to undertake such an investigation in the near
future, fully expecting to see all of the non-axisymmetric dynamics,
including the emergence of vortex breakdown of S-type at the trailing edge
of the bubble - which may signal transition to turbulent flow regimes, and
Shilnikov chaos.

Finally, there are those who might simply say that the ideal fluid flow
context (enabling us to reap the considerable benefits of a Hamiltonian
structure for the kinematics) of our B-type breakdown model disqualifies it
from serious consideration of such phenomena in real fluids. Certainly there
is considerable truth to this with regard to some flow quantities; however,
extensive experimental and numerical studies of B-type vortex breakdown in
both real, high Reynolds number flows and ideal flows have shown that both
B-type and S-type vortex breakdown behavior tends to be qualitatively the
same in all of these cases, and even quite close quantitatively speaking
over time intervals of several minutes. Thus the ideal tells us a great deal
- especially dynamically - about the real when it comes to vortex breakdown.

\section*{Acknowledgments}

The main part of this paper was completed while D. Blackmore was a visiting
professor at the Technical University of Denmark (DTU), with the M. Br\o ns
serving as his immediate host and principal scientific collaborator. D.
Blackmore wants to especially thank his hosts at DTU for providing him with
such a stimulating and intellectually rich environment, along with a chance
to enjoy the extraordinary Danish hospitality, and a wonderful opportunity
to reconnect with his roots. He also wishes to acknowledge the generous
support from the The Villum Kann Rasmussen Foundation during his visit to
DTU. Debts of gratitude are also owed by D. Blackmore to Egon Krause for
many insightful conversations and much encouragement that were instrumental
in the formulation of the vortex breakdown model, and to Lu Ting and Omar
Knio for sharing their considerable expertise and wisdom regarding the work
in this paper.

\end{document}